\documentclass[a4paper,11pt]{amsart}
\usepackage{amssymb}
\usepackage{amsmath}
\usepackage{amsthm}
\usepackage[pdftex]{graphicx}
\usepackage[english]{babel}
\usepackage{array}
\usepackage{titletoc}
\usepackage[T1]{fontenc}  
\usepackage[applemac]{inputenc} 
\usepackage[all]{xy}
\usepackage[vcentermath]{youngtab}
\usepackage{stmaryrd}

\usepackage[mathcal]{euscript}

 \numberwithin{equation}{section}

  \theoremstyle{plain}
\newtheorem{theoreme}[equation]{Theorem}
\newtheorem{lemme}[equation]{Lemma}
\newtheorem{lemme-def}[equation]{Lemma-Definition}
\newtheorem{proposition}[equation]{Proposition}
\newtheorem{corollaire}[equation]{Corollary}
  \theoremstyle{definition}
\newtheorem{definition}[equation]{Definition}

  \theoremstyle{remark}
\newtheorem{remarque}[equation]{Remark}
\newtheorem{notations}[equation]{Notations}
\newtheorem{exemple}[equation]{Example}

 \newcommand{\theo}{\begin{theoreme}}
 \newcommand{\defi}{\begin{definition}}
 \newcommand{\rema}{\begin{remarque}}
 \newcommand{\prop}{\begin{proposition}}
 \newcommand{\coro}{\begin{corollaire}}
 \newcommand{\lemm}{\begin{lemme}}
 \newcommand{\exem}{\begin{exemple}}
 \newcommand{\nota}{\begin{notations}}

 \newcommand{\etheo}{\end{theoreme}}
 \newcommand{\edefi}{\end{definition}}
 \newcommand{\erema}{\end{remarque}}
 \newcommand{\eprop}{\end{proposition}}
 \newcommand{\ecoro}{\end{corollaire}}
 \newcommand{\elemm}{\end{lemme}}
 \newcommand{\eexem}{\end{exemple}}
 \newcommand{\enota}{\end{notations}}

 \newcommand{\becen}{\begin{center}}
 \newcommand{\ecen}{\end{center}}

 \newcommand{\benu}{\begin{enumerate}}
 \newcommand{\eenu}{\end{enumerate}}

 \newcommand{\bite}{\begin{itemize}}
 \newcommand{\eite}{\end{itemize}}

\def\[#1\]{\begin{align*}#1\end{align*}}

\newcommand{\demo}{\begin{proof}}
\newcommand{\edemo}{\end{proof}}

\newcommand{\CC}{\mathbb C}
\newcommand{\QQ}{\mathbb Q}
\newcommand{\ZZ}{\mathbb Z}
\newcommand{\NN}{\mathbb N}

\newcommand{\Hom}{\operatorname{Hom}}
\newcommand{\End}{\operatorname{End}}

\newcommand{\aaa}{\alpha}
\newcommand{\bbb}{\beta}
\newcommand{\ccc}{\gamma}
\newcommand{\eee}{\epsilon}
\newcommand{\Spec}{\operatorname{Spec}}

\newcommand{\Irr}{\operatorname{Irr}}

\newcommand{\rarrow}{\rightarrow}

\newcommand{\Rarrow}{\Rightarrow}
\newcommand{\Lrarrow}{\Leftrightarrow}
\newcommand{\codim}{\operatorname{codim}}

 \newcommand{\dsum}{\displaystyle\sum}
 \newcommand{\doplus}{\displaystyle\bigoplus}

 \newcommand{\dcup}{\displaystyle\bigcup}
 \newcommand{\dcap}{\displaystyle\bigcap}

\newcommand{\isom}{\overset{\sim}{\rightarrow}}

\newcommand{\id}{\operatorname{id}}
\newcommand{\ima}{\operatorname{Im}}
\newcommand{\tra}{\operatorname{Tr}}
\newcommand{\W}{\textup{\sffamily W}}
\newcommand{\M}{\textup{\sffamily M}}

\newcommand{\La}{\textup{\sffamily L}}
\newcommand{\Z}{\textup{\sffamily Z}}

\author{Tristan Bozec} 

\title{Quivers with loops and generalized crystals}
\date{}

\begin{document}

\maketitle 

\begin{abstract} 
In the context of varieties of representations of arbitrary quivers, possibly carrying loops, we define a generalization of Lusztig Lagrangian subvarieties. From the combinatorial study of their irreducible components arises a structure richer than the usual Kashiwara crystals. Along with the geometric study of Nakajima quiver varieties, in the same context, this yields a notion of generalized crystals, coming with a tensor product. As an application, we define the semicanonical basis of the Hopf algebra generalizing quantum groups, which was already equipped with a canonical basis. The irreducible components of the Nakajima varieties provide the family of highest weight crystals associated to dominant weights, as in the classical case.
\end{abstract}

\tableofcontents

\section*{Introduction}
Lusztig defined in~\cite{lulu} Lagrangian subvarieties of the cotangent stack to the moduli stack of representations of a quiver associated to any Kac-Moody algebra. The proof of the Lagrangian character of these varieties was obtained via the study of some natural stratifications of each irreducible component, and then proceeding by induction. The particular combinatorial structure thus attached to the set of irreducible components made it possible for Kashiwara and Saito in~\cite{kashisaito} to relate this variety to the usual quantum group associated to Kac-Moody algebras, via the notion of \emph{crystals}. This later led Lusztig in~\cite{semicanonical} to define a \emph{semicanonical basis} of this quantum group, indexed by the irreducible components of these Lagrangian varieties.

There are more and more evidence of the relevance of the study of quivers with loops. A particular class of such quivers are the comet-shaped quivers, which have recently been used by Hausel, Letellier and Rodriguez Villegas in their study of the topology of character varieties, where the number of loops at the central vertex is the genus of the considered curve (see~\cite{MR2453601} and~\cite{MR3003926}). We can also see quivers with loops appearing in a work of Nakajima relating quiver varieties with branching (see~\cite{MR2470410}), as in the work of Okounkov and Maulik about quantum cohomology (see~\cite{OM}).

Kang, Kashiwara and Schiffmann generalized these varieties in the framework of generalized Kac-Moody algebra in~\cite{MR2553376}, using quivers with loops. In this case, one has to impose a somewhat unnatural restriction on the regularity of the maps associated to the loops.

In this article we define a generalization of such Lagrangian varieties in the case of arbitrary quivers, possibly carrying loops. As opposed to the Lagrangian varieties constructed by Lusztig, which consisted in nilpotent representations, we have to consider here slightly more general representations. That this is necessary is already clear from the Jordan quiver case. Note that our Lagrangian variety is strictly larger than the one considered in~\cite{MR2553376} and has many more irreducible components. 
Our proof of the Lagrangian character is also based on induction, but with non trivial first steps, consisting  in the study of quivers with one vertex but possible loops. From our proof emerges a new combinatorial structure on the set of irreducible components, which is more general than the usual crystals, in that there are now more operators associated to a vertex with loops, see~\ref{generA}.

In a second section we study Nakajima varieties, still in the context of arbitrary quivers. We construct Lagrangian subvarieties, and generalize the notion of tensor product of their irreducible components, introduced by Nakajima in~\cite{nakatens}. The geometric statements obtained in the two first sections give the intuition of the way crystals and their tensor product should be generalized, which is done in a third part. The algebraic definition and study of the crystal $\mathcal B(\infty)$ enable us to define a semicanonical basis for the positive part of the generalized quantum group $U^+$ defined in~\cite{article2}, where it is already equipped with a canonical basis, built \textit{via} the theory of Lusztig perverse sheaves associated to quivers with loops.
We finally use our study of Nakajima quiver varieties to produce a geometric realization of the generalized crystals $\mathcal B(\lambda)$.

\subsubsection*{Acknowledgement}
I would like to thank Olivier Schiffmann for his constant support and availability during the preparation of this work.

\section{Lusztig quiver varieties}

Let $Q$ be a quiver, defined by a set of vertices $I$ and a set of oriented edges $\Omega=\{h:s(h)\rarrow t(h)\}$. We denote by $\bar h:t(h)\rarrow s(h)$ the opposite arrow of $h\in\Omega$, and $\bar Q$ the quiver $(I,H=\Omega\sqcup\bar\Omega)$, where $\bar\Omega=\{\bar h\mid h\in\Omega\}$: each arrow is replaced by a pair of arrows, one in each direction, and we set  $\eee(h)=1$ if $h\in \Omega$, $\eee(h)=-1$ if $h\in \bar\Omega$. Note that the definition of $\bar h$ still makes sense if $h\in\bar\Omega$.
We denote by $\Omega(i)$ the set of loops of $\Omega$ at $i$, and call $i$ \emph{imaginary} if $\omega_i=|\Omega(i)|≥1$, \emph{real} otherwise. Denote by $I^\textup{im}$ (resp. $I^\textup{re}$) the set of imaginary vertices (resp. real vertices).
Finally, set $H(i)=\Omega(i)\sqcup\bar\Omega(i)$.

We work over the field of complex numbers $\CC$.

For any pair of $I$-graded $\CC$-vector spaces $V=(V_i)_{i\in I}$ and $V'=(V'_i)_{i\in I}$, we set:\[
\bar E(V,V')=& \bigoplus_{h\in  H}\Hom(V_{s(h)},V'_{t(h)}).\]
For any dimension vector $\nu=(\nu_i)_{i\in I}$, we fix an $I$-graded $\CC$-vector space $V_\nu$ of dimension $\nu$, and put $\bar E_\nu=\bar E(V_\nu,V_\nu)$.
The space $\bar E_\nu=\bar E(V_\nu,V_\nu)$ is endowed with a symplectic form:\[
\omega_\nu(x,x')=\dsum_{h\in H}\tra(\epsilon(h)x_hx'_{\bar h})\]
which is preserved by the natural action of $G_\nu=\prod_{i\in I} GL_{\nu_i}(\CC)$ on $\bar E_\nu$. The associated moment map $\mu_\nu:\bar E_\nu\rarrow\mathfrak g_\nu= \oplus_{i\in I}\End(V_\nu)_i$ is given by:\[
\mu_\nu(x)=\dsum_{h\in H}\epsilon(h)x_{\bar h}x_h.\]
Here we have identified $\mathfrak g_\nu^*$ with $\mathfrak g_\nu$ via the trace pairing.

\defi\label{nillu}
An element $x \in \bar E_\nu$ is said to be \emph{seminilpotent} if there exists an $I$-graded flag $\W=(\W_0=V_\nu\supset\ldots\supset\W_r=\{0\})$ of $V_\nu$ such that:\[
\left.\begin{array}{ll}  
x_h(\W_\bullet)\subseteq \W_{\bullet+1}&\text{if }h\in \Omega,\\
x_h(\W_\bullet)\subseteq \W_{\bullet}&\text{if }h\in \bar\Omega.\end{array}\right.\]
We put $\Lambda(\nu)=\{x\in \mu_\nu^{-1}(0)\mid x\text{ seminilpotent}\}$. 
\edefi

\lemm\label{isotropie} The variety $\Lambda(\nu)$ is isotropic.\elemm

\demo
We use the following general fact (see \textit{e.g.}~\cite[§8.3]{MR1299726}): 
\prop Let $X$ be a smooth algebraic variety, $Y$ a projective variety and $Z$ a smooth closed algebraic subvariety of $X\times Y$. Consider the Lagrangian subvariety $\Lambda=T^*_Z(X\times Y)$ of $T^*(X\times Y)$. Then the image of the projection $q:\Lambda\cap( T^*X\times T^*_YY)\rarrow T^*X$ is isotropic.\eprop

We apply this result to $X=\oplus_{h\in \Omega}\End(V_{\nu_{s(h)}},V_{\nu_{t(h)}})$, $Y$ the $I$-graded flag variety of $V_\nu$ and:\[
Z=\{(x,\W)\in X\times Y\mid x(\W_{\bullet})\subseteq\W_{\bullet+1}\} .\]
In this case, we get:\[
T^*X&=\bar E_\nu  \\
T^*Y&=\{(\W,\xi)\in Y\times \mathfrak g_\nu\mid \xi(\W_\bullet)\subseteq\W_{\bullet+1}\}  \\
\Lambda&=
\left\{
(x,\W,\xi)
~\middle|~
\begin{aligned} 
&\xi=\dsum_{h\in H}\eee(h)x_{\bar h}x_h  \\ 
&\forall h\in \Omega,~ x_h(\W_\bullet)\subseteq \W_{\bullet+1}\text{ and }x_{\bar h}(\W_\bullet)\subseteq \W_{\bullet}
\end{aligned}
\right\}\\
\ima q&=\left\{x\in \bar E_\nu
~\middle|~\begin{aligned}& \mu_\nu(x)=0\text{ and there exists }\W\in Y \text{ such that} \\
&\forall h\in \Omega,~ x_h(\W_\bullet)\subseteq \W_{\bullet+1}\text{ and }x_{\bar h}(\W_\bullet)\subseteq \W_{\bullet}\end{aligned}\right\},
\]
hence $\Lambda(\nu)\subseteq\ima q$, which proves the lemma.
\edemo

\subsection{The case of the Jordan quiver}
This case is very well known. For $\nu\in\NN$, we have:\[
\Lambda(\nu)=\{(x,y)\in(\End\CC^\nu)^2\mid x\text{ nilpotent and}~[x,y]=0\}=\dcup_{\lambda}T^*_{\mathcal O_\lambda}(\End\CC^\nu),\]
where $\mathcal O_\lambda$ is the nilpotent orbit associated to the partition $\lambda$ of $\nu$. Therefore $\Lambda(\nu)$ is a Lagrangian subvariety of $(\End\CC^\nu)^2$, and its irreducible components are the closures of the conormal bundles to the nilpotent orbits.

\subsection{The case of the quiver with one vertex and $g≥2$ loops}
For $\nu\in\NN$, $\Lambda(\nu)$ is the subvariety of $(\End\CC^\nu)^{2g}$ with elements $(x_i,y_i)_{1≤i≤g}$ such that:\bite
\item[$\triangleright$] there exists a flag $\W$ of $\CC^\nu$ such that $x_i(\W_\bullet)\subseteq \W_{\bullet+1}$ and $y_i(\W_\bullet)\subseteq \W_\bullet$;
\item[$\triangleright$] $\dsum_{1≤i≤g}[x_i,y_i]=0$.\eite

We will denote by $\mathcal C_\nu$ the set of compositions of $\nu$, \textit{i.e.}~tuples $\textup{\sffamily c}=(\textup{\sffamily c}_1,\ldots,\textup{\sffamily c}_r)$ of $\NN_{>0}$ such that:\[
|\textup{\sffamily c}|=\sum_{1≤k≤r}\textup{\sffamily c}_k=\nu.\]

We will also often forget the index $1≤i≤g$ in the rest of this section, which is dedicated to the proof of the following theorem:

\theo\label{multithm} The subvariety $\Lambda(\nu)\subseteq(\End\CC^\nu)^{2g}$ is Lagrangian. Its irreducible components are parametrized by $\mathcal C_\nu$. 
\etheo

\nota

For $(x_i,y_i)\in\Lambda(\nu)$, we define $\W_0(x_i,y_i)=\CC^\nu$, then by induction $\W_{k+1}(x_i,y_i)$ the smallest subspace of $\CC^\nu$ containing $\sum x_i(\W_k(x_i,y_i))$ and stable by $(x_i,y_i)$. By seminilpotency, we can define $r$ to be the first power such that $\W_r(x_i,y_i)=\{0\}$. Although $r$ depends on $(x_i,y_i)$ we don't write it explicitly. 

Let $\textup{\sffamily c}(x_i,y_i)$ denotes the composition associated to the flag $\W_\bullet(x_i,y_i)$:\[
\textup{\sffamily c}_k(x_i,y_i)=\dim\dfrac{\W_{k-1}(x_i,y_i)}{\W_{k}(x_i,y_i)}.\]
For every $\textup{\sffamily c}\in\mathcal C_\nu$, we define a locally closed subvariety: \[
\Lambda(\textup{\sffamily c})=\left\{(x_i,y_i)\in\Lambda(\nu)\mid \dim\dfrac{\W_{\bullet-1}(x_i,y_i)}{\W_{\bullet}(x_i,y_i)}=\textup{\sffamily c}\right\}\subseteq\Lambda(\nu).\]
 Then, if $\delta=(\delta_1,\ldots,\delta_{r-1})\in\NN^{r-1}$, let $\Lambda(\textup{\sffamily c})_{\delta}\subseteq\Lambda(\textup{\sffamily c})$ be the locally closed subvariety defined by:\[
\dim\bigg(\dcap_{1≤i≤g}\ker
\left\{\xi  \mapsto   y_i^{(k)}\xi-\xi y_i^{(k+1)}\right\}\bigg)=\delta_k
,\]
where:\[
y_i^{(k)}\in\End\left(\dfrac{\W_{k-1}(x_i,y_i)}{\W_{k}(x_i,y_i)}\right)\]
 is induced by $y_i$ and:\[
\xi\in\Hom\left(\dfrac{\W_{k}(x_i,y_i)}{\W_{k+1}(x_i,y_i)},\dfrac{\W_{k-1}(x_i,y_i)}{\W_{k}(x_i,y_i)}\right).\]
Set $l=\textup{\sffamily c}_1$, then:\[
\check\Lambda(\textup{\sffamily c})_{\delta}=\left\{(x_i,y_i,\mathfrak X,\bbb,\ccc)~\middle|~\begin{aligned}&(x_i,y_i)\in\Lambda(\textup{\sffamily c})_{\delta}\\
&\W_{1}(x_i,y_i)\oplus\mathfrak X=\CC^\nu
\\&\bbb:\W_{1}(x_i,y_i)\isom\CC^{\nu-l}\text{ and }  \ccc:\mathfrak X\isom\CC^{l} \end{aligned}\right\},\]
and:\[
\pi_{\textup{\sffamily c},\delta}\left|\begin{aligned}& \check\Lambda(\textup{\sffamily c})_{\delta}\rarrow{\Lambda(\textup{\sffamily c}^-)_{\delta^-}}\times(\End\CC^{l})^g
 \\  &(x_i,y_i,\mathfrak X,\bbb,\ccc)\mapsto(\bbb_*(x_i,y_i)_{\W_{1}},\ccc_*({y_i})_{\mathfrak X}))
  \end{aligned}\right.\]
where $\textup{\sffamily c}^-=(\textup{\sffamily c}_2,\ldots,\textup{\sffamily c}_{r})$ and $\delta^-=(\delta_2,\ldots,\delta_{r-1})$. Let finally $(\Lambda(\textup{\sffamily c}^-)_{\delta^-}\times(\End\CC^{l})^g)_{\textup{\sffamily c},\delta}$ denote the image of $\pi_{\textup{\sffamily c},\delta}$.
\enota

\prop\label{irreduc} The morphism $\pi_{\textup{\sffamily c},\delta}$ is smooth over its image, with connected fibers of dimension $\nu^2+(2g-1)l(\nu-l)+\delta_{1}$ whenever $\Lambda(\textup{\sffamily c})_\delta\neq\varnothing$.\eprop

\demo
Let $(x_i,y_i,z_i)\in({\Lambda(\textup{\sffamily c}^-)_{\delta^-}}\times(\End\CC^{l})^g)_{\textup{\sffamily c},\delta}$. Let $\mathfrak W$ and $\mathfrak X$ be two supplementary subspaces of $\CC^\nu$ such that $\dim\mathfrak X=l$, together with two isomorphisms: \[
\bbb:\mathfrak W\isom\CC^{\nu-l}\text{ and }\ccc:\mathfrak X\isom\CC^{l}
.\]
We identify $x_i,y_i$ and $z_i$ with $\bbb^*(x_i,y_i)$ and $\ccc^*z_i$, and define an element $(X_i,Y_i)$ in the fiber of $(x_i,y_i,z_i)$ by setting:\[
(X_i,Y_i)_{\mathfrak W}&=(x_i,y_i)\\
(X_i,Y_i)_{\mathfrak X}&=(0,z_i)\\
(X_i,Y_i)_{|\mathfrak X}^{|\mathfrak W}&=(u_i,v_i)\in\Hom(\mathfrak X,\mathfrak W)^{2g}.
\]
Then:\[
\mu_\nu(X_i,Y_i)=0\Lrarrow\phi(u_i,v_i)=\dsum_{i=1}^g(x_iv_i+u_iz_i-y_iu_i)=0,\]
and, for $\xi \in\Hom(\mathfrak W,\mathfrak X )$:\[
 \forall(u_i,v_i),~\tra(\xi \phi(u_i,v_i))=0 & \Lrarrow   \left\{\begin{aligned}& \forall i,\forall u_i,~\tra(\xi (u_iz_i-y_iu_i))=0 \\  &\forall i,\forall v_i,~\tra(\xi x_iv_i)=0\end{aligned}\right.  \\
  &  \Lrarrow  \left\{\begin{aligned}& \forall i,\forall u_i,~\tra((z_i\xi -\xi y_i)u_i)=0 \\  & \forall i,\forall v_i,~\tra(\xi x_iv_i)=0\end{aligned}\right.  \\
    & \Lrarrow   \left\{\begin{aligned}&\forall i,~ z_i\xi =\xi y_i \\  & \forall i,~\xi x_i=0\end{aligned}\right. \\
    & \Lrarrow
   \left\{\begin{aligned}& \W_{1}(x_i,y_i)\subseteq\ker \xi  \\ &  \forall i,~ z_i\xi ^{(1)}=\xi ^{(1)}y_i^{(1)} \end{aligned}\right.  
\]
where $\xi ^{(1)}$ denotes the map $\mathfrak W/\W_{1}(x_i,y_i)\rarrow\mathfrak X$ induced by $\xi $. Since $(x_i,y_i,z_i)$ is in the image of $\pi_{\textup{\sffamily c},\delta}$, the image of $\phi$ is of codimension $\delta_{1}$, and thus its kernel is of dimension $(2g-1)l(\nu-l)+\delta_{1}$.

Moreover, if we denote by $u_i^{(1)}$ the map $\mathfrak X\rarrow\mathfrak W/\W_{1}(x_i,y_i)$ induced by $u_i$,  $\W_1(X_i,Y_i)=\mathfrak W$ if and only if the space spanned by the action of $(y_i^{(1)})_i$ on $\sum_i \ima u_i^{(1)}$ is $\mathfrak W/\W_{1}(x_i,y_i)$. This condition defines an open subset of $\ker\phi$.

We end the proof noticing that the set of elements $(\mathfrak W,\mathfrak X,\bbb,\ccc)$ is isomorphic to $GL_\nu(\CC)$.
\edemo

\prop\label{nonvacuite} The variety $\Lambda(\textup{\sffamily c})_0$ is not empty.\eprop

\demo Fix $\W$ of dimension $\textup{\sffamily c}$ and define $x_1$ such that \[
x_1(\W_\bullet)&\subseteq\W_{\bullet+1}\\
 {x_1}_{|\W_{k-1}/\W_k}^{|\W_{k}/\W_{k+1}}&\neq0.\] 
 We define inductively an element $y_1$ stabilizing $\W$ such that:\benu
 \item[$\triangleright$] the action of ${y_1}^{(k+1)}$ on $\ima\left( {x_1}_{|\W_{k-1}/\W_k}^{|\W_{k}/\W_{k+1}}\right)$ spans $\W_{k}/\W_{k+1}$;
 \item[$\triangleright$] $\Spec{y_1}^{(k)}\cap\Spec{y_1}^{(k+1)}=\varnothing$.
 \eenu
 We finally set $x_2=-x_1,y_2=y_1$ and $x_i=y_i=0$ for $i>2$. This yields an element $(x_i,y_i)$ in $\Lambda(\textup{\sffamily c})_0$.
\edemo

\coro\label{dimirred} For any $\textup{\sffamily c}\in\mathcal C_\nu$, $\Lambda(\textup{\sffamily c})_0$ is irreducible of dimension $g\nu^2$.\ecoro

\demo We argue by induction on $r$. If $\textup{\sffamily c}=(\nu)$, we have $\Lambda(\textup{\sffamily c})_0=\Lambda(\textup{\sffamily c})=(\End\CC^\nu)^g$ which is irreducible of dimension $g\nu^2$. For the induction step,~\ref{irreduc} and~\ref{nonvacuite} ensure us that $\check{\Lambda}(\textup{\sffamily c})_0$ is irreducible of dimension:\begin{multline*}
\nu^2+(2g-1)l(\nu-l)+\dim(\Lambda(\textup{\sffamily c}^-)_{0}\times(\End\CC^{l})^g)_{\textup{\sffamily c},0}\\=\nu^2+(2g-1)l(\nu-l)+g(\nu-l)^2+gl^2\end{multline*}
since $(\Lambda(\textup{\sffamily c}^-)_{0}\times(\End\CC^{l})^g)_{\textup{\sffamily c},0}$ is a non empty subvariety of $\Lambda(\textup{\sffamily c}^-)_{0}\times(\End\CC^{l})^g$, irreducible of dimension $g(\nu-l)^2+gl^2$ by our induction hypothesis. Moreover,\[
\check{\Lambda}(\textup{\sffamily c})_0\rarrow{\Lambda(\textup{\sffamily c})_0}\]
being a principal bundle with fibers of dimension $\nu^2-l(\nu-l)$, we get that ${\Lambda(\textup{\sffamily c})_0}$ is irreducible of dimension \[
\nu^2+(2g-1)l(\nu-l)+g(\nu-l)^2+gl^2-\nu^2+l(\nu-l)=g\nu^2.\]
\edemo 

\lemm Let $V$ and $W$ be two vector spaces, and $k≥0$. For any $(u,v)\in\End V\times\End W$, we set:\[
\mathcal C(u,v)&=\{x\in\Hom(V,W)\mid xu=vx\}\\
(\End V\times\End W)_k&=\{(u,v)\in\End V\times\End W\mid \dim\mathcal C(u,v)=k\}.\]
Then we have \[\codim (\End V\times\End W)_k≥k.\]
\elemm

\demo
The restriction of an endomorphism $a$ to a generalized eigenspace associated to an eigenvalue $\eta$ will be denoted by $a_\eta=\eta\id+\tilde a_\eta$. As usual, the nilpotent orbit associated to a partition $\xi$ will be denoted by $\mathcal O_\xi$.
 We have:\[
&\codim (\End V\times\End W)_k\\
&=\codim\{(u,v)\mid \textstyle\sum_{\aaa,\bbb}\dim\mathcal C(u_\aaa,v_\bbb)=k\}\\
&=\codim\{(u,v)\mid \textstyle\sum_{\aaa\in\Spec u\cap\Spec v}\dim\mathcal C(u_\aaa,v_\aaa)=k\}\\
&=\codim\{(u,v)\mid \textstyle\sum_{\aaa}\dim\mathcal C(\tilde u_\aaa,\tilde v_\aaa)=k\}\\
&=\codim\left\{(u,v)\left|\begin{array}{l} (\tilde u_\aaa,\tilde v_\aaa)\in\mathcal O_{\lambda_\aaa}\times\mathcal O_{\mu_\aaa} \\ \sum_\aaa\sum_j(\lambda'_\aaa)_j(\mu_\aaa')_j=k\end{array}\right.\right\}
\]
Thus, \[
&\codim (\End V\times\End W)_k≥k\\
&\Lrarrow \sum_\aaa(\codim\mathcal O_{\lambda_\aaa}+\codim\mathcal O_{\mu_\aaa}-1)≥\sum_\aaa\textstyle\sum_j(\lambda'_\aaa)_j(\mu'_\aaa)_j\\
&\Lrarrow\sum_\aaa(\textstyle{\sum_j(\lambda'_\aaa)_j^2+\sum_j(\mu'_\aaa)^2_j}-1)≥\displaystyle\sum_\aaa\textstyle\sum_j(\lambda'_\aaa)_j(\mu'_\aaa)_j,\]
which is clear.
\edemo

\prop If $\delta\neq0$, we have $\dim\Lambda(\textup{\sffamily c})_\delta<g\nu^2$.\eprop

\demo It's enough to show that if $\delta_{1}>0$, we have:\[
\dim({\Lambda(\textup{\sffamily c}^-)_{\delta^-}}\times(\End\CC^{l})^g)_{\textup{\sffamily c},\delta}+\delta_{1}<\dim({\Lambda(\textup{\sffamily c}^-)_{0}}\times(\End\CC^{l})^g).\]
This is a consequence of the previous lemma (recall that $g≥2$). Indeed, if we set :\[
((\End V)^g\times(\End W)^g)_k=\{(u_i,v_i)\mid \dim\cap_i\mathcal C(u_i,v_i)=k\},\]
we have:\[
((\End V)^g\times(\End W)^g)_k\subseteq\prod_{i=1}^g(\End V\times\End W)_{k_i}\]
for some $k_i≥k$, and thus:\begin{multline*}
\codim((\End V)^g\times(\End W)^g)_k\\≥\sum_i\codim(\End V\times\End W)_{k_i}≥\sum_ik_i≥gk>k.\end{multline*}
\edemo

The following proposition concludes the proof of theorem~\ref{multithm}:

\prop Every irreducible component of $\Lambda(\textup{\sffamily c})$ is of dimension larger than or equal to $g\nu^2$.\eprop

\demo We first prove the result for the following variety:\[
\tilde\Lambda(\textup{\sffamily c})=\{((x_i,y_i),\W)\in\Lambda(\nu)\times Y_\textup{\sffamily c}\mid x_i(\W_\bullet)\subseteq\W_{\bullet+1}\text{ and }y_i(\W_\bullet)\subseteq\W_{\bullet}\}\]
where $Y_\textup{\sffamily c}$ denotes the variety of flags of $\CC^\nu$ of dimension w.
We use the following notations, analogous to~\ref{isotropie}:\[
X&=\{(x_i)_{1≤i≤g}\in(\End\CC^\nu)^g\}\\
Z&=\{((x_i)_{1≤i≤g},\W)\mid x_i(\W_\bullet)\subseteq\W_{\bullet+1}\}\subseteq X\times Y_\textup{\sffamily c}.\]
We get:\[
T^*X&=\{(x_i,y_i)_{1≤i≤g}\in(\End\CC^\nu)^{2g}\}\\
T^*Y_\textup{\sffamily c}&=\{(\W,K)\in Y_\textup{\sffamily c}\times\End\CC^\nu\mid K(\W_\bullet)\subseteq\W_{\bullet+1}\}\]
and:\[
T^*_Z(X\times Y_\textup{\sffamily c})&=\left\{((x_i,y_i),\text{F},K)~\middle|~\begin{aligned}& \sum_{1≤i≤g}[x_i,y_i]=K  \\ &
x_i(\W_\bullet)\subseteq\W_{\bullet+1}\text{ and }y_i(\W_\bullet)\subseteq\W_{\bullet}\end{aligned}\right\}\]
which is a pure Lagrangian subvariety of $T^*(X\times Y_\textup{\sffamily c})$, of dimension $g\nu^2+\dim Y_\textup{\sffamily c}$. Since $T^*Y_\textup{\sffamily c}$ is irreducible of dimension $2\dim Y_\textup{\sffamily c}$, the irreducible components of the fibers of $T^*_Z(X\times Y_\textup{\sffamily c})\rarrow T^*Y_\textup{\sffamily c}$ are of dimension larger than or equal to $g\nu^2-\dim Y_\textup{\sffamily c}$. We denote by $\tilde\Lambda_{\text{\sffamily W}}$ the fiber above $(\W,0)$, and by $P$ the stabilizer of $\W$ in $G_\nu$.
Since $G_\nu$ and $P$ are irreducible, we get that the components of:\[
\tilde\Lambda(\textup{\sffamily c})=G_\nu{\times}_P\tilde\Lambda_{\text{\sffamily W}}\]
are of dimension larger than or equal to $\dim Y_\textup{\sffamily c}+(g\nu^2-\dim Y_\textup{\sffamily c})=g\nu^2$.

We extend this result to $\Lambda(\textup{\sffamily c})$, noticing that:\[
\left.\begin{array}{rcl} \Lambda(\textup{\sffamily c}) & \hookrightarrow  & \tilde\Lambda(\textup{\sffamily c})   \\
(x_i,y_i)  &  \mapsto & (x_i,y_i,\W_\bullet(x_i,y_i)) \end{array}\right.\]
defines an open embedding.
\edemo

\subsection{The general case}

Denote by $a_{i,j}$ the number of edges of $\Omega$ such that $s(h)=i$ and $t(h)=j$, and denote by:\[ 
C=(2\delta_{i,j}-a_{i,j}-a_{j,i})\]
 the Cartan matrix of $Q$.
For every $\nu,\bbb\in\NN^I$ and $j\in I$, we put:\begin{align*}
\langle\nu,\bbb\rangle&=\dsum_{i\in I}\nu_i\bbb_i\\
e_j&=(\delta_{i,j})_{i\in I}.\end{align*}

\defi For every subset $i\in I$, and every $x\in\Lambda(\nu)$, we denote by $
\mathfrak I_i(x)$ the subspace of $V_\nu$ spanned by the action of $x$ on $\oplus_{j\neq i}V_j$. Then, for $l>0$, we set:\[
\Lambda(\nu)_{i, l}=\{x\in\Lambda(\nu)\mid\codim\mathfrak I_i(x)= le_i\}.\]
\edefi

\rema\label{remlu} By the definition of seminilpotency, we have:\[
\Lambda(\nu)=\dcup_{i\in I,l≥1}\Lambda(\nu)_{i,l}.\]
Indeed, if $x\in\Lambda(\nu)$, there exists an $I$-graded flag $(\W_0\supset\ldots\supset\W_r)$ such that $(x,\W)$ satisfies~\ref{nillu}. Therefore there exists $i\in I$ and $l>0$ such that $\W_0/\W_{1}\simeq V_{le_i}$, and thus $x\in\cup_{k≥l}\Lambda(\nu)_{i,k}$.
\erema

\prop\label{generA} There exists a variety $\check\Lambda(\nu)_{i, l}$ and a diagram:
\[\xymatrix{
&\check{\Lambda}(\nu)_{i, l}\ar[ld]_{q_{i, l}}\ar[rd]^{p_{i, l}}&\\
\Lambda(\nu)_{i, l}&&\Lambda(\nu-le_i)_{i, 0}\times\Lambda(le_i)}\]
such that $p_{i, l}$ and $q_{i, l}$ are smooth with connected fibers,  inducing a bijection:\[
\Irr\Lambda(\nu)_{i, l}{~\isom~}\Irr\Lambda(\nu-le_i)_{i, 0}\times\Irr\Lambda(le_i).\]
\eprop

\demo In this proof we will denote by $I(V,V')$ the set of $I$-graded isomorphisms between two $I$-graded spaces $V$ and $V'$ of same $I$-graded dimension.
We set:\[
\check\Lambda(\nu)_{i, l}=\left\{(x,\mathfrak X,\bbb,\ccc)~\middle|~\begin{aligned}& x\in\Lambda(\nu)_{i, l} \\  & \mathfrak X \text{ $I$-graded and } \mathfrak I_i(x)\oplus\mathfrak X=V_\nu\\  &
\bbb\in I(\mathfrak I_i(x), V_{\nu-le_i} )\text{ and }\ccc\in I(\mathfrak X, V_{le_i}) \end{aligned}\right\}\]
and:\[
p_{i, l}\left|\begin{aligned}& \check\Lambda(\nu)_{i, l}   \rarrow  \Lambda(\nu-le_i)_{i, 0}\times\Lambda(le_i)  \\  &   (x,\mathfrak X,\bbb,\ccc)  \mapsto  
(\bbb_*(x_{\mathfrak I_i(x)}),\ccc_*(x_{\mathfrak X})). \end{aligned}\right.\]
We study the fibers of $p_{i, l}$: take $y\in\Lambda(\nu-le_i)_{i, 0}$ and $z\in\Lambda(le_i)$ and consider $\mathfrak I$ and $\mathfrak X$ two supplementary $I$-graded subspaces of $V_\nu$ such that $\dim\mathfrak X=le_i$, together with two isomorphisms:\[
\bbb\in I(\mathfrak I, V_{\nu-le_i} )\text{ and }\ccc\in I(\mathfrak X, V_{le_i}).\]
 We identify $y$ and $z$ with $\bbb^*y$ and $\ccc^*z$, and we define a preimage $x$ by setting $x_{|\mathfrak I}^{|\mathfrak I}=y$, $x_{|\mathfrak X}^{|\mathfrak X}=z$ and $x_{|\mathfrak X}^{|\mathfrak I}=\eta\in \bar E(\mathfrak X,\mathfrak I)$. In order to get $\mu_\nu(x)=0$, $\eta$ must satisfy the following relation:\[
\phi(\eta)=\dsum_{h\in H:s(h)=i}\epsilon( h)(y_{\bar h}\eta_h+\eta_{\bar h}z_h)=0.\]
We need to show that $\phi$ is surjective to conclude. Consider $\xi\in \Hom(\mathfrak I_i,\mathfrak X_i)$ such that $\operatorname{Tr}(\phi(\eta)\xi)=0$ for every $\eta$.
For every edge $h$ such that $s(h)=i\neq j=t(h)$ and every $\eta_h$, we have:\[
0&=\operatorname{Tr}(y_{\bar h}\eta_h\xi)\\
&=\operatorname{Tr}(\xi y_{\bar h}\eta_h)\]
Hence $\xi y_{\bar h}=0$, and $\ima y_{\bar h}\subseteq\ker\xi$.
Now consider a loop $h\in H(i)$. For every $\eta_h$, we have:\[
0&=\operatorname{Tr}\big((\eta_hz_{\bar h}-y_{\bar h}\eta_h)\xi\big)\\
&=\operatorname{Tr}\big(\eta_h(z_{\bar h}\xi-\xi y_{\bar h})\big).\]
Hence $\xi y_{\bar h}=z_{\bar h}\xi $ and therefore $\ker \xi$ is stable by $y_{\bar h}$. As $\codim\mathfrak I_i(y)=0$, we get $\xi=0$, which finishes the proof.
\edemo

 We can now state the following theorem, which answers a question asked in~\cite{yli}:
 
\theo The subvariety $\Lambda(\nu)$ of $\bar E_\nu$ is Lagrangian.\etheo

\demo
Since this subvariety is isotropic by~\ref{isotropie} we just have to show that the irreducible components of $\Lambda(\nu)$ are of dimension $\langle\nu,(1-C/2)\nu\rangle$. We proceed by induction on $\nu$, the first step corresponding to the one vertex quiver case which has already been treated: we have seen that $\Lambda({le_i})$ is of dimension $\langle le_i,(1-C/2)le_i\rangle$.

Next, consider $C\in\Irr\Lambda(\nu)$ for some $\nu$. By~\ref{remlu}, there exists $i\in I$ and $l≥1$ such that $C\cap\Lambda(\nu)_{i,l}$ is dense in $C$. Let $\check C=(C_1,C_2)$ the couple of irreducible components corresponding to $C$ via the bijection obtained in~\ref{generA}:\[
\Irr \Lambda(\nu)_{i,l}\isom\Irr \Lambda(\nu-le_i)_{i,0}\times\Irr\Lambda(le_i) .\]
We also know by the proof of~\ref{generA} that the fibers of $p_{i, l}$ are of dimension:\[
\langle \nu,\nu\rangle+\langle \nu-le_i,(1-C)le_i\rangle.
\]
Since $q_{i, l}$ is a principal bundle with fibers of dimension $\langle\nu,\nu\rangle-\langle le_i,\nu-le_i\rangle$, we get: \[
\dim C=\dim\check C+\langle \nu-le_i,(2-C)le_i\rangle.\]
But $\Lambda(\nu-le_i)_{i,0}$ is open in $\Lambda(\nu-le_i)$, so we can use our induction hypothesis and the first step to write: \[\dim\check C=\langle \nu-le_i,(1-C/2)(\nu-le_i)\rangle+\langle le_i,(1-C/2)le_i\rangle\]
and thus obtain:\[
\dim C=\langle \nu,(1-C/2)\nu\rangle.\]
\edemo

\subsection{Constructible functions}

Following~\cite{semicanonical}, we denote by $\mathcal M(\nu)$ the $\QQ$-vector space of constructible functions $\Lambda(\nu)\rarrow\QQ$, which are constant on any $G_\nu$-orbit. Put $\mathcal M=\oplus_{\nu≥0}\mathcal M(\nu)$, which is a graded algebra once equipped with the product $*$ defined in~\cite[2.1]{semicanonical}.

For $Z\in \Irr\Lambda(\nu)$ and $f\in\mathcal M(\nu)$, we put $\rho_Z(f)=c$ if $Z\cap f^{-1}(c)$ is an open dense subset of $Z$. 

If $i\in I^\textup{im}$ and $(l)$ denotes the trivial composition or partition of $l$, we denote by $1_{i,l}$ the characteristic function of the associated irreducible component $Z_{i,(l)}\in\Irr\Lambda(le_i)$ (the component of elements $x$ such that $x_h=0$ for any loop $h\in\Omega(i)$). If $i\notin I^\textup{im}$, we just denote by $1_i$ the function mapping to $1$ the only point in $\Lambda(e_i)$.

We have $1_{i,l}\in\mathcal M(le_i)$ for $i\in I^\textup{im}$ and $1_i\in\mathcal M(e_i)$ for $i\notin  I^\textup{im}$. We denote by $\mathcal M_\circ\subseteq\mathcal M$ the subalgebra generated by these functions.

\lemm\label{lemmconstr} Suppose $Q$ has one vertex $\circ$ and $g≥1$ loop(s). For every $Z\in \Irr\Lambda(\nu)$ there exists $f\in\mathcal M_\circ(\nu)$ such that $\rho_Z(f)=1$ and $\rho_{Z'}(f)=0$ for $Z'\neq Z$.\elemm

\demo We denote by $Z_\textup{\sffamily c}$ the irreducible component associated to the partition (resp. composition) $\textup{\sffamily c}$ of $\nu$ if $g=1$ (resp. $g≥2$). By convention, if $g=1$, $Z_\textup{\sffamily c}$ will denote the component associated to the orbit $\mathcal O_\textup{\sffamily c}$ defined by:\[
x\in\mathcal O_\textup{\sffamily c}\Lrarrow\dim\ker x^i=\sum_{1≤k≤i}\textup{\sffamily c}_k.\]

If $g≥2$, we remark that by trace duality, we can assume that $Z_\textup{\sffamily c}$ is the closure of $\Lambda\check{}_\textup{\sffamily c}$ defined by:\[
(x_i,y_i)_{1≤i≤g}\in \Lambda\check{}_\textup{\sffamily c}\Lrarrow\dim\text K_i=\sum_{1≤k≤i}\textup{\sffamily c}_k\]
where we define by induction K$_0=\{0\}$, then K$_{j+1}$ as the biggest subspace of  $\cap_ix_i^{-1}(\text K_j)$ stable by $(x_i,y_i)$.
From now on, $\textup{\sffamily c}=(\textup{\sffamily c}_1,\ldots,\textup{\sffamily c}_r)$ will denote indistinctly a partition or a composition depending on the value of $g$. 
We define an order by:\[
\textup{\sffamily c}\preceq{\textup{\sffamily c}'}\text{ if and only if for any $i≥1$ we have }\sum_{1≤k≤i}\textup{\sffamily c}_k≤\sum_{1≤k≤i}{\textup{\sffamily c}'}_k.\]
Therefore, setting $\tilde 1_\textup{\sffamily c}=1_{\textup{\sffamily c}_r}*\cdots*1_{\textup{\sffamily c}_1}$, where $1_{l}=1_{\circ,l}$, we get:\[ 
x\in Z_\textup{\sffamily c},~\tilde1_{\textup{\sffamily c}'}(x)\neq0~\Rarrow~{\textup{\sffamily c}'}\preceq\textup{\sffamily c}.\]
For $\textup{\sffamily c}=(\nu)$ we have $\tilde1_\textup{\sffamily c}=1_{\nu}$ which is the characteristic function of $Z_\textup{\sffamily c}$, and we put $1_\textup{\sffamily c}=\tilde1_\textup{\sffamily c}$ in this case. Then, by induction:\[ 
1_{\textup{\sffamily c}}=\tilde 1_\textup{\sffamily c}-\dsum_{{\textup{\sffamily c}'}\prec\textup{\sffamily c}}\rho_{Z_{\textup{\sffamily c}'}}(\tilde 1_\textup{\sffamily c})1_{{\textup{\sffamily c}'}}\]
has the expected property.
\edemo

\nota {~}
\bite
\item[$\triangleright$] From now on, if $\textup{\sffamily c}$ corresponds to an irreducible component of $\Lambda(|\textup{\sffamily c}|e_i)$, we will note $1_{i,\textup{\sffamily c}}$ the function corresponding to $1_{\textup{\sffamily c}}$ in the previous proof.
\item[$\triangleright$] For $Z\in \Irr\Lambda(\nu)_{i,l}$, we denote by $\eee_i(Z)\in\Irr\Lambda(le_i)$ the composition of the second projection with the bijection obtained in~\ref{generA}. Note that $|\eee_i(Z)|=l$.\eite\enota

\prop For every $Z\in \Irr\Lambda(\nu)$, there exists $f\in\mathcal M_\circ(\nu)$ such that $\rho_Z(f)=1$ and $\rho_{Z'}(f)=0$ if $Z'\neq Z$.\eprop

\demo We proceed as in~\cite[lemma 2.4]{semicanonical}, by induction on $\nu$. The first step consists in~\ref{lemmconstr}. Then, consider $Z\in \Irr\Lambda(\nu)$. There exists $i\in I$ and $l>0$ such that $Z\cap\Lambda(\nu)_{i,l}$ is dense in $Z$. 

We now proceed by descending induction on $l$. There's nothing to say if $l>\nu_i$.

Otherwise, let $(Z',Z_\textup{\sffamily c})\in\Irr\Lambda(\nu-le_i)_{i,0}\times\Irr\Lambda(le_i)$ be the pair of components corresponding to $Z$. By the induction hyopthesis on $\nu$, there exists $g\in\mathcal M_\circ(\nu-le_i)$ such that $\rho_{\overline {Z'}}(g)=1$ and $\rho_{Y}(g)=0$ if $\overline{Z'}\neq Y\in\Irr\Lambda(\nu-le_i)$.

Then we set $\tilde f=1_{i,\textup{\sffamily c}}*g\in\mathcal M_\circ(\nu)$, and get:\bite
\item[$\bullet$] $\rho_Z(\tilde f)=1$,
\item[$\bullet$] $\rho_{Z'}(\tilde f)=0$ if $Z'\in\Irr\Lambda(\nu)\setminus Z$ satisfies $|\eee_i(Z')|=l$,
\item[$\bullet$] $\tilde f(x)=0$ if $x\in\Lambda(\nu)_{i,<l}$ so that $\rho_{Z'}(\tilde f)=0$ if $|\eee_i(Z')|<l$. \eite
If $|\eee_i(Z')|>l$, we use the induction hypothesis on $l$: there exists $f_{Z'}\in\mathcal M_\circ(\nu)$ such that $\rho_{Z'}(f_{Z'})=1$ and $\rho_{Z''}(f_{Z'})=0$ if $Z''\in\Irr\Lambda(\nu)\backslash Z'$. We end the proof by setting:\[
f=\tilde f-\dsum_{Z': |\eee_i(Z')|>l}\rho_{Z'}(\tilde f)f_{Z'}.\]
\edemo

\section{Nakajima quiver varieties}

Fix an $I$-graded vector space $W$ of dimension $\lambda=(\lambda_i)_{i\in I}$.
For any dimension vector $\nu=(\nu_i)_{i\in I}$, we still fix an $I$-graded $\CC$-vector space $V_\nu=((V_\nu)_i=V_{\nu_ie_i})_{i\in I}$ of dimension $\nu$.
We will denote by $(x,f,g)=((x_h)_{h\in H},(f_i)_{i\in I},(g_i)_{i\in I})$ the elements of the following space:\[
E(V,\lambda)= \bar E(V,V)\oplus\bigoplus_{i\in I}\Hom(V_i,W_i)\bigoplus_{i\in I}\Hom(W_i,V_i)\]
defined for any $I$-graded space $V$, and put $E_{\nu,\lambda}=E(V_\nu,\lambda)$ for any dimension vector $\nu$.
This space is endowed with a symplectic form:\[
\omega_{\nu,\lambda}\big((x,f,g),(x',f',g')\big)=\dsum_{h\in H}\tra(\epsilon(h)x_hx'_{\bar h})+\dsum_{i\in I}\tra(g_if'_i-g'_if_i)\]
which is preserved by the natural action of $G_\nu=\prod_{i\in I} GL_{\nu_i}(\CC)$ on $E_{\nu,\lambda}$. The associated moment map $\mu_{\nu,\lambda}:E_{\nu,\lambda}\rarrow\mathfrak g_\nu= \oplus_{i\in I}\End(V_\nu)_i$ is given by:\[
\mu_{\nu,\lambda}(x,f,g)=\Big(g_if_i+\dsum_{h\in H: s(h)=i}\epsilon(h)x_{\bar h}x_h\Big)_{i\in I}.\]
Here we have identified $\mathfrak g_\nu^*$ with $\mathfrak g_\nu$ via the trace pairing.
Put:\[
\M_\circ(\nu,\lambda)=\mu_{\nu,\lambda}^{-1}(0).\]

\defi Set $\chi:G_\nu\rarrow\CC^*$, $(g_i)_{i\in I}\mapsto\prod_{i\in I}\det^{-1}g_i$. We denote by:\[
\mathfrak M_\circ(\nu,\lambda)&=\M_\circ(\nu,\lambda)/\!\!/ G_\nu\\
\mathfrak M(\nu,\lambda)&=\M_\circ(\nu,\lambda)/_{\!\chi} G_\nu\]
the geometric and symplectic quotients (with respect to $\chi$).\edefi

\prop An element $(x,f,g)\in \M_\circ(\nu,\lambda)$ is stable with respect to $\chi$ if and only if the only $x$-stable subspace of $\ker f$ is $\{0\}$. Set:\[
\M(\nu,\lambda)=\{(x,f,g)\in \M_\circ(\nu,\lambda)\mid(x,f,g)\textup{ stable}\},\]
then $\mathfrak M(\nu,\lambda)=\M(\nu,\lambda)/\!\!/ G_\nu$.\eprop

\subsection{A crystal-type structure}

\defi\label{nillunak}
An element $(x,f,g) \in E_{\nu,\lambda}$ is said to be \emph{seminilpotent} if $x\in\bar E_\nu$ is, according to~\ref{nillu}.
We put:\[
\La_\circ(\nu,\lambda)=\{(x,f,0)\in \M_\circ(\nu,\lambda)\mid x\textup{ seminilpotent}\}\subseteq \M_\circ(\nu,\lambda)\]
 and define $\La(\nu,\lambda)\subseteq \M(\nu,\lambda)$ in the same way. Finally set:\[
 \mathfrak L_\circ(\nu,\lambda)&=\La_\circ(\nu,\lambda)/\!\!/G_\nu\\
 \mathfrak L(\nu,\lambda)&=\La_\circ(\nu,\lambda)/_{\!\chi} G_\nu=\La(\nu,\lambda)/\!\!/G_\nu.\]
 We will simply denote by $(x,f)$ the elements of $\La_\circ(\nu,\lambda)$. 
\edefi

 There is an alternative definition of $\mathfrak L(\nu,\lambda)$. Define a $\CC^*$-action on $\mathfrak M(\nu,\lambda)$ by:\[
t\diamond[x,f,g]=[t^{(1+\eee)/2}x,f,tg].\]
We have:\[
\mathfrak L(\nu,\lambda)=\{[x,f,g]\mid\exists\lim_{t\to\infty}t\diamond[x,f,g]\}.\]
By the same arguments than in~\cite[5.8]{instan}, we have the following:
\prop \label{nakalag}
The subvariety ${\mathfrak L}(\nu,\lambda)\subset{\mathfrak M}(\nu,\lambda)$ is Lagrangian.
\eprop
Note that since we consider seminilpotents instead of nilpotents, we still have:\[
\omega_{\nu,\lambda}(t\diamond-,-)=t\omega_{\nu,\lambda}.\]

\defi For every subset $i\in I$, and every $(x,f,g)\in \M_\circ(\nu,\lambda)$, we denote by $
\mathfrak I_i(x,f,g)$ the subspace of $V_{\nu}$ spanned by the action of $x\oplus g$ on $(\oplus_{j\neq i}V_j)\oplus W_i$. Then, for $l≥0$, we set:\[
\M_\circ(\nu,\lambda)_{i, l}=\{x\in \M_\circ(\nu,\lambda)\mid\codim\mathfrak I_i(x,f,g)= le_i\}.\]
We define $\M(\nu,\lambda)_{i,l}$, $\La_\circ(\nu,\lambda)_{i,l}$ and $\La(\nu,\lambda)_{i,l}$ in the same way. The quantity $\codim\mathfrak I_i(x,f,g)$ being stable on $G_\nu$-orbits, the notations $\mathfrak M_\circ(\nu,\lambda)_{i,l}$, $\mathfrak M(\nu,\lambda)_{i,l}$, $\mathfrak L_\circ(\nu,\lambda)_{i,l}$ and $\mathfrak L(\nu,\lambda)_{i,l}$ also make sense.
\edefi

\rema\label{remnak} ~

\bite\item As in~\ref{remlu}, we have:\[
 \La_\circ(\nu,\lambda)=\displaystyle\bigsqcup_{i\in I,l≥1} \La_\circ(\nu,\lambda)_{i,l}.\]
\item Note that $\La_\circ(le_i,0)=\Lambda({le_i})$.\eite\erema

\prop\label{fibnak} There exists a variety $\check\M_\circ(\nu,\lambda)_{i,l}$ and a diagram:
\begin{align}\label{fibnakdiag}\xymatrix{& ~~~\check{\M}_\circ(\nu,\lambda)_{i,l}\ar[ld]_{q_{i,l}}\ar[rd]^{ p_{i,l}}&\\\M_\circ(\nu,\lambda)_{i,l}&&\M_\circ(\nu-le_i,\lambda)_{i,0}\times \M_\circ({le_i},0)
}\end{align}
such that $p_{i, l}$ and $q_{i, l}$ are smooth with connected fibers,  inducing a bijection:\[
\Irr\M_\circ(\nu,\lambda)_{i,l}{~\isom~}\Irr\M_\circ(\nu-le_i,\lambda)_{i,0}\times\Irr\M_\circ({le_i},0).\]
\eprop

\demo In this proof we will denote by $I(V,V')$ the set of $I$-graded isomorphisms between two $I$-graded spaces $V$ and $V'$ of same $I$-graded dimension.
We set:\[
\check\M_\circ(\nu,\lambda)_{i,l}=\left\{(x,f,g,\mathfrak X,\bbb,\ccc)~\middle|~\begin{aligned}& (x,f,g)\in\M_\circ(\nu,\lambda)_{i,l} \\  & \mathfrak X\textup{ $I$-graded and }\mathfrak I_i(x,f,g)\oplus\mathfrak X=V_\nu\\  &
\bbb\in I(\mathfrak I_i(x,f,g), V_{\nu-le_i} )\\ &\ccc\in I(\mathfrak X, V_{le_i}) \end{aligned}\right\}\]
and:\[
p_{i, l}\left|\begin{aligned}& \check\M_\circ(\nu,\lambda)_{i,l}  \rarrow  \M_\circ(\nu-le_i,\lambda)_{i,0}\times \M_\circ({le_i},0)  \\  &   (x,f,g,\mathfrak X,\bbb,\ccc)  \mapsto  
(\bbb_*(xf,g)_{\mathfrak I_i(x,f,g)},\ccc_*(x,f,g)_{\mathfrak X}). \end{aligned}\right.\]
We study the fibers of $p_{i,l}$: take $(x,f,g)\in\M_\circ(\nu-le_i,\lambda)_{i,0}$ and $(z,0,0)\in\M_\circ({le_i},0)$ and consider $\mathfrak I$ and $\mathfrak X$ two supplementary $I$-graded subspaces of $V_\nu$ such that $\dim\mathfrak X=le_i$, together with two isomorphisms:\[
\bbb\in I(\mathfrak I, V_{\nu-le_i} )\textup{ and }\ccc\in I(\mathfrak X, V_{le_i}).\]
 We identify $(x,f,g)$ and $z$ with $\bbb^*(x,f,g)$ and $\ccc^*z$, and we define a preimage $(X,F,G)$ by setting $(X,F,G)_{|\mathfrak I\oplus W}^{|\mathfrak I\oplus W}=(x,f,g)$, $X_{|\mathfrak X}^{|\mathfrak X}=z$ and:\[
 (X,F)_{|\mathfrak X}^{|\mathfrak I\oplus W}=(\eta,\theta)
 \in \bar E(\mathfrak X,\mathfrak I)\oplus\Hom(\mathfrak X_i,W_i).\]
 In order to get $\mu_{\nu,\lambda}(X,F,G)=0$, $(\eta,\theta)$ must satisfy the following relation:\[
\psi(\eta,\theta)=\dsum_{h\in H:s(h)=i}\epsilon( h)(x_{\bar h}\eta_h+\eta_{\bar h}z_h)+g_i\theta_i=0.\]
We need to show that $\psi$ is surjective to conclude. Consider $\xi\in \Hom(\mathfrak I_i,\mathfrak X_i)$ such that $\operatorname{Tr}(\psi(\eta,\theta)\xi)=0$ for every $(\eta,\theta)$. Then we have for every edge $h\in H$ such that $s(h)=i\neq j=t(h)$ and for every $\eta_h$:\[
0&= \operatorname{Tr}(x_{\bar h}\eta_h\xi)\\
&=\operatorname{Tr}(\eta_h\xi x_{\bar h}).\]
 Hence $\xi x_{\bar h}=0$ and $\ima x_{\bar h}\subseteq\ker\xi$. We also have $\operatorname{Tr}(g_i\theta_i\xi)=0$ for every $\theta_i$, so we similarly get $\ima g_i\subseteq\ker\xi$. Now consider a loop $h\in H$ at $i$. We have for every $\eta_h$:\[
0&= \operatorname{Tr}\big((x_{\bar h}\eta_h-\eta_hz_{\bar h})\xi\big)\\
&=\operatorname{Tr}\big(\eta_h(\xi x_{\bar h}-z_{\bar h}\xi)\big),\]
hence $\xi x_{\bar h}=z_{\bar h}\xi$ and therefore $\ker \xi$ is stable by $x_{\bar h}$. 
 Since $(x,f,g)\in\M_\circ(\nu-le_i,\lambda)_{i,0}$, we get $\xi=0$, which finishes the proof.
\edemo

	\coro We also have a bijection:\[
\textup{\sffamily l}_\circ(\nu,\lambda)_{i,l}:	\Irr\La_\circ(\nu,\lambda)_{i,l}{~\isom~}\Irr\La_\circ(\nu-le_i,\lambda)_{i,0}\times\Irr\La_\circ({le_i},0).\]
\ecoro

\demo The image of a seminilpotent element by $p_{i,l}$ is a pair of seminilpotent elements, and the fiber of $p_{i,l}$ over a pair of seminilpotent elements consists in seminilpotent elements.\edemo

\subsection{Extension to the stable locus}

We will often use the following well-known fact:
\lemm\label{lemm0}  Consider $y\in\End \mathfrak I$ and $z\in\End\mathfrak X$ such that $\Spec y\cap\Spec z=\varnothing$. If $\CC[y].v=\mathfrak I$ and $\CC[z].v'=\mathfrak X$ 
for some $v\in\mathfrak I$ and $v'\in\mathfrak X$, then $
\CC[ y\oplus z].v\oplus v'=\mathfrak I\oplus\mathfrak X$.\elemm

\nota Let $i$ be imaginary and put $\Omega(i)=\{b_1,\ldots,b_{\omega_i}\}$. For every $(x,f)\in\La_\circ(\nu,\lambda)$, we set $\sigma_i(x)=x^*_{\bar b_1}$, where $*$ stands for the duality:\[
\End V&\rarrow\End V^*=\End\big(\Hom(V,\CC)\big)\\
u&\mapsto u^*=[\phi\mapsto\phi\circ u]\]
for every $\CC$-vector space $V$.
\enota

\lemm\label{lemm1} For every $C\in\Irr\Lambda(le_i)$, there exists $x\in C$ such that:\[
\exists \psi\in V^*_{le_i}, \CC[ \sigma_i(x)].\psi=V^*_{le_i}.\]
\elemm

\demo It's a consequence of sections 1.1 and 1.2. If $\omega_i=1$ and $\lambda$ is a partition of $l$, denote by $\mu$ the conjugate partition of $\lambda$. Let $x\in\mathcal O_\lambda$ be defined in a base:\[
e=(e_{1,1},\ldots,e_{1,\mu_1},\ldots,e_{r,1},\ldots,e_{r,\mu_r})\]
by:\[
x^*_{b_1}= 
  \left(\!\!\!\!\!\!
     \raisebox{0.5\depth}{
       \xymatrixcolsep{1ex}
       \xymatrixrowsep{1ex}
       \xymatrix{
         J_{\mu_1} \ar @{.}[dddrrr] & 
         0 \ar @{.}[rr]\ar @{.}[ddrr]&
         &  0 \ar@{.}[dd]\\
         0 \ar@{.}[dd] \ar@{.}[ddrr]&\\
         &&& 0\\
         0 \ar@{.}[rr]  & & 0 & J_{\mu_r}
       }
     }
   \right)
   \textup{ and }
x^*_{\bar b_1}= 
  \left(\!\!\!\!\!\!
     \raisebox{0.5\depth}{
       \xymatrixcolsep{1ex}
       \xymatrixrowsep{1ex}
       \xymatrix{
         t_1I_{\mu_1}+J_{\mu_1} \ar @{.}[dddrrr] & 
         0 \ar @{.}[rr]\ar @{.}[ddrr]&
         &  0 \ar@{.}[dd]\\
         0 \ar@{.}[dd] \ar@{.}[ddrr]&\\
         &&& 0\\
         0 \ar@{.}[rr]  & & 0 & t_rI_{\mu_r}+J_{\mu_r}
       }
     }
   \right)
  \]
where the $t_i$ are all distinct and nonzero, and:\[
J_p=\left(\!\!\!\!\!\!
     \raisebox{0.5\depth}{
       \xymatrixcolsep{1ex}
       \xymatrixrowsep{1ex}
       \xymatrix{
         0 \ar @{.}[dddrrr] & 
         1 \ar @{.}[ddrr]
         & 0 \ar @{.}[r] \ar@{.}[dr]&  0 \ar@{.}[d]\\
         0 \ar@{.}[dd] \ar@{.}[ddrr]&&&0\\
         &&& 1 \\
         0 \ar@{.}[rr] &  & 0 & 0
       }
     }
   \right).\]
It is enough to consider $\psi$ with nonzero coordinates relatively to $(e_{1,\mu_1},\ldots,e_{r,\mu_r})$ to get $\CC[\sigma_i(x)].\psi=V^*_{le_i}$.
   If $\omega_i≥2$, we use the proof of~\ref{nonvacuite}: in any irreducible component we can define $x$ such that there exists $v$ such tat $\CC[ x_{\bar b_1}].v=V_{le_i}$ ($x_{\bar b_i}$ corresponds to $y_i$ in the aforementioned proof, $x_{ b_i}$ to $x_i$). We get the result by duality.
\edemo

\rema Note that the case $\omega_i=1$ is very well known since it corresponds to the case of the Hilbert scheme of points in the plane.\erema

\defi Set:\[
\La(\lambda):=\dcup_{\nu}\La(\nu,\lambda)\subseteq\dcup_\nu\La_\circ(\nu,\lambda)=:\La_\circ(\lambda),\]
and define $\textup B(\lambda)$ as the smallest subset of $\Irr\La_\circ(\lambda)$ containing the only element of $\Irr\La_\circ(0,\lambda)$, and stable by the $\textup{\sffamily l}_\circ(\nu,\lambda)^{-1}_{i,l}(-,\Irr\Lambda(le_i))$ for $\nu,i,l$ such that:\bite
\item $\langle e_i,\lambda-C\nu\rangle≥-l$ if $i\in I^\textup{re}$,
\item $\lambda_i+\sum_{h\in H_i} \nu_{t(h)}>0$ if $i\in I^\textup{im}$
\eite
where $H_i=\{h\in H\mid i=s(h)\neq t(h)\}$.
\edefi

\lemm\label{lemmclef} For every $i\in I^\textup{im}$, we write $\Omega(i)=\{b_{i,1},\ldots, b_{i,\omega_i}\}$. For every $C\in \textup B(\lambda)$, there exists $(x,f)\in C$ such that:\begin{align}\label{specialcondition}
\left\{\begin{aligned} & (x,f)\textup{ stable} \\
  & \forall i\in I^\textup{im},\exists \phi_i\in W^*_i\oplus(\oplus_{h\in H_i}V^*_{\nu_{t(h)}}), \CC[\sigma_i(x)].\Sigma_i(x,f)(\phi_i)=V^*_{\nu_i} \end{aligned}\right.
\end{align}
where $\Sigma_i(x,f)=f^*_i+\sum_{h\in H_i}x^*_h$.
\elemm

\demo
We proceed by induction on $\nu$, with first step consisting in the case of $C\in \textup B(\lambda)\cap \Irr\La_\circ(le_i,\lambda)$ for some $l>0$. If $i\notin I^\textup{im}$, we have $l≤\lambda_i$ by definition of $\textup B(\lambda)$, hence we can find $(x,f)\in C$ such that~\ref{specialcondition} since it's equivalent here to $f$ injective. If $i\in I^\textup{im}$, we have $\lambda_i>0$ by definition of $\textup B(\lambda)$, and we can use~\ref{lemm1}.

Now consider $C\in\textup B(\lambda)\cap\Irr\La_\circ(\nu,\lambda)_{i,l}$ for some $\nu$ and $l>0$, and set $(C_1,C_2)=\textup{\sffamily l}_\circ(\nu,\lambda)_{i,l}(C)$. First assume that $i\notin I^\textup{im}$. Thanks to the induction hypothesis, we can pick $((x,f),z)\in C_1\times C_2$ such that $(x,f)$ satisfies~\ref{specialcondition}. Following the notations used in the proof of~\ref{fibnak}, we build an element of $C$ satisfying~\ref{specialcondition} by chosing $(\eta,\theta)$ such that $\theta+\sum_{h\in H_i}\eta_h$ is injective with values in a supplementary of $\ima(f_i+\sum_{h\in H_i}x_h)$ in $W_i\oplus\ker(\sum_{h\in H_i}x_{\bar h})$: it's possible since $l+\langle e_i,\lambda-C\nu\rangle≥0$ by definition of $\textup B(\lambda)$.

If $i\in I^\textup{im}$, take $(x,f)\in C_1$ satisfying~\ref{specialcondition} and $z\in C_2$ such that:\[
\left\{\begin{aligned}& \Spec x_{\bar b_{i,1}}\cap\Spec z_{\bar b_{i,1}}=\varnothing \\  
&\exists \psi\in V^*_{le_i},\CC[\sigma_i(z)].\psi=V^*_{le_i},  \end{aligned}\right.
\]
which is possible, thanks to~\ref{lemm1}.
Still following the notations of the proof of~\ref{fibnak}, we build an element of $C$ mapped to $((x,f),z)$ by considering $(\eta,\theta)$ such that:\[
\Big(\theta^*+\sum_{h\in H_i}\eta^*_h\Big)(\phi_i)=\psi\]
where $\phi_i\in W^*_i\oplus(\oplus_{h\in H_i}V^*_{\nu_{t(h)}})$ satisfies $\CC[\sigma_i(x)].\Sigma_i(x,f)(\phi_i)=\mathfrak I^*$ (we use the induction hypothesis), which is possible even if $\mathfrak I=\{0\}$ since we have $W^*_i\oplus(\oplus_{h\in H_i}V^*_{\nu_{t(h)}})\neq\{0\}$ by definition of $\textup B(\lambda)$. Put $\eta_{b_{i,j}}=\eta_{\bar b_{i,j}}=0$ for every $j≥2$, so that:\[
\psi_i(\eta,\theta)=0\Lrarrow x_{\bar b_{i,1}}\eta_{b_{i,1}}-\eta_{b_{i,1}}z_{\bar b_{i,1}}=\dsum_{h\in H_i}\eee(h)(x_{\bar h}\eta_h+\eta_{\bar h}z_h).\]
Hence we can choose $\eta_{b_{i,1}}$ in order to satisfy the right hand side equation since:\[
\Spec x_{\bar b_{i,1}}\cap\Spec z_{\bar b_{i,1}}=\varnothing\Rarrow (\eta_{b_{i,1}}\mapsto x_{\bar b_{i,1}}\eta_{b_{i,1}}-\eta_{b_{i,1}} z_{\bar b_{i,1}})\textup{ invertible}.\]
Thanks to~\ref{lemm0}, $(X,F)\in C$ satisfies:\[
\CC[\sigma_i(X)].\Sigma_i(X,F)(\phi_i)=V^*_{\nu_i}.\]
 We finally have to check the stability of $(X,F)$ to conclude. Consider $S\subseteq\ker F$ stable by $X$. We have $S\cap{\mathfrak I}=\{0\}$ by stability of $(x,f)$, thus $S\simeq S_i$ and we see $S$ as a subspace of $\ker F\cap(\cap_{h\in H_i}\ker X_h)$. But then $S^*$ is stable by $\sigma_i(X)$ and contains $\ima F^*+\sum_{h\in H_i}\ima X^*_h$, and thus $\phi_i$. Hence $S^*=V_{\nu_i}$, and $S=\{0\}$.
\edemo

\prop\label{subcryst} We have $\textup B(\lambda)=\Irr\La(\lambda)$.
\eprop

\demo Thanks to~\ref{lemmclef}, we have $\textup B(\lambda)\subseteq\Irr\La(\lambda)$. Consider $Z\in\Irr\La(\nu,\lambda)_{i,l}\setminus \textup B(\lambda)$ for some $l>0$. We know (\textit{c.f.}~\cite[4.6]{nakajima}) that if $i\in I^\textup{re}$, we necessarily have $l+\langle e_i,\nu-C\lambda\rangle≥0$, and thus, by definition of $\textup B(\lambda)$:\[
\textup{\sffamily l}_\circ(\nu,\lambda)_{i,l}(Z)\in\Big(\Irr\La(\nu-le_i,\lambda)\setminus \textup B(\lambda)\Big)\times\Irr\Lambda(le_i).\]
 If $i\in I^\textup{im}$, $Z\in\Irr\La(\nu,\lambda)_{i,l}$ necessarily implies $\lambda_i+\sum_{h\in H_i} \nu_{t(h)}>0$, and we get to the same conclusion. By descending induction on $\nu$, we obtain that the only irreducible component of $\La(0,\lambda)$ doesn't belong to $\textup B(\lambda)$, which is absurd.
\edemo

\coro\label{bijquot} Take $i\in I^\textup{im}$ and assume $\Irr\La(\nu,\lambda)_{i,l}\subseteq \textup B(\lambda)$. We have the following commutative diagram:
\begin{align}\label{bijquotdiag}\xymatrix{
 \Irr\La(\nu,\lambda)_{i,l}\ar[rr]_-\sim^-{\textup{\sffamily l}(\nu,\lambda)_{i,l}}\ar[d]_{\sim}&& \Irr{\La}({\nu-le_i}, \lambda)_{i,0}\times \Irr\Lambda(le_i)\ar[d]^{\sim}\\
  \Irr{\mathfrak L}(\nu,\lambda)_{i,l}\ar[rr]_-{\sim}^-{\mathfrak l(\nu,\lambda)_{i,l}}&& \Irr{\mathfrak L}({\nu-le_i}, \lambda)_{i,0}\times \Irr\Lambda(le_i).
}\end{align}
\ecoro

\demo
By definition of stability, the action of $G_{\nu}$ on $\La(\nu,\lambda)$ is free.
\edemo

\subsection{Tensor product on $\Irr\mathfrak L$}

\subsection{Another Lagrangian subvariety}

Embed $W$ in a $\lambda+\lambda'$-dimensional $I$-graded vector space, and fix a supplementary subspace $W'$ of $W$. We still denote by $I(X,Y)$ the set of $I$-graded isomorphisms between two $I$-graded spaces $X$ and $Y$.

For every $\mathbf v\in\NN^I$, denote by $\Z_\circ(\mathbf v)\subseteq\M_\circ(\mathbf v,\lambda+\lambda')$ of elements $(x,f,g)$ such that there exists an $I$-graded subspace $V$ of $V_\mathbf v$ satisfying:\benu
\item $x(V)\subseteq V$;
\item $f(V)\subseteq W$;
\item $g(W\oplus W')\subseteq V$;
\item $g(W)=\{0\}$,\eenu
and denote by $V(x,f,g)$ the larger $x$-stable subspace of $f^{-1}(W)$ containing $\ima g$.
We will then denote by $\widetilde\Z_\circ(\mathbf v)\subset\Z_\circ(\mathbf v)$ the subvariety of elements $(x,f,g)$ such that:\[
({x},{f})_{| V\times  V}^{| V\times W}\textup{ and }({x},{f})_{|(V_{\mathbf v}/ V)\times(V_{\mathbf v}/ V)}^{|(V_{\mathbf v}/ V)\times(W\oplus W'/W)}\textup{ are seminilpotents}\]
where we have written $V$ instead of $V(x,f,g)$.
We get a stratification of $\widetilde\Z_\circ(\mathbf v)$ by setting, for any $\nu,\nu'$ such that $\nu+\nu'=\mathbf v$:\[
\widetilde\Z_\circ(\nu,\nu')=\left\{(x,f,g)\in\widetilde\Z_\circ(\nu+\nu')\mid\dim V(x,f,g)=\nu\right\}.\]
Define the following incidence variety:\[ 
\check\Z_\circ(\nu,\nu')=
\left\{(x,f,g,V',\bbb)~\middle|~\begin{aligned}&(x,f,g)\in\widetilde\Z_\circ(\nu,\nu')\\ & V(x,f,g)\oplus V'=V_{\nu+\nu'}\\  &
 \bbb\in I( V(x,f,g),V_\nu)\times I(V',V_{\nu'})
\end{aligned}\right\}.\]
By definition of $ V(x,f,g)$ (again denoted by $ V$ hereunder), we have:\[
(x,f,g)\in\Z_\circ(\mathbf v)\Rarrow ({x},{f})_{|(V_{\mathbf v}/ V)\times(V_{\mathbf v}/ V)}^{|(V_{\mathbf v}/ V)\times(W\oplus W'/W)}\textup{ stable},\]
hence the following application is well defined:\[
\textup{\sffamily T}_\circ\left|\begin{aligned}&\check\Z_\circ(\nu,\nu')  \rarrow    \La_\circ(\nu,\lambda)\times\La(\nu',\lambda') \\  &  (x,f,g, V',\bbb)  \mapsto   \bbb_*\left(({x},{f})_{| V\times  V}^{| V\times W},({x},{f})_{| V'\times V'}^{| V'\times(W\oplus W'/W)}\right) \end{aligned}\right.
\]

\prop\label{fibtens} The map $\textup{\sffamily T}_\circ$ is smooth with connected fibers.\eprop

\demo
Let $(x,f)$ and $(x',f')$ be elements of $\La_\circ(\nu,\lambda)$ and $\La(\nu',\lambda')$ and take $I$-graded spaces $V$ and $V'$ of dimensions $\nu$ and $\nu'$. Define $(X,F,G,V',\bbb)$ in the fiber $\textup{\sffamily T}_\circ^{-1}((x,f),(x',f'))$ by:
\benu\itemsep0.5em
\item $\bbb\in I(V,V_\nu)\times I(V',V_{\nu'})$;
\item $G=0\oplus\tau$ where:\[
\nu\in \oplus_{i\in I}\Hom(W'_i,V_i);\]
\item $X=\bbb^*x\oplus(\bbb^*x'+\eta):V\oplus V'\rarrow V\oplus V'$ where:\[\eta\in \oplus_{h\in H}\Hom(V'_{s(h)},V_{t(h)});\]
\item $F=\bbb^*f\oplus(\bbb^*f'+\theta):V\oplus V'\rarrow W\oplus W'$ where:\[
 \theta\in \oplus_{i\in I}\Hom(V'_i,W_i) ;\]
\eenu
such that $\mu_{\nu+\nu',\lambda+\lambda'}(X,F,G)=0$. 

\lemm\label{lemmfibtens} This equation is  linear in the variables $(\tau,\eta,\theta)$, and the associated linear map is surjective.\elemm

\demo We first identify $x,x'$, and $f'$ with $\bbb^*x,\bbb^*x'$, and $\bbb^*f'$. Then the linear map $\zeta=(\zeta_i)$ we are interested in is given by:\[
\zeta_i(\tau,\eta,\theta)=\tau_if'_i+\dsum_{h\in H:s(h)=i}\epsilon(\bar h)(x_{\bar h}\eta_h+\eta_{\bar h}x'_h).\]
Take $L\in \oplus_{i\in I}\Hom(V_i,V'_i)$ such that for every $(\tau,\eta,\theta)$:\[
\sum_{i\in I}\operatorname{Tr}(\zeta(\tau,\eta,\theta)L_i)=0.\]
Then for every edge $h$ such that $s(h)=i$, $t(h)=j$, we have for every $\eta_h$:\[
\operatorname{Tr}(x_{\bar h}\eta_hL_i)- \operatorname{Tr}(\eta_hx'_{\bar h}L_j)=0.\]
But\[
\operatorname{Tr}(\eta_hL_ix_{\bar h})- \operatorname{Tr}(\eta_hx'_{\bar h}L_j)=\operatorname{Tr}(\eta_hL_ix_{\bar h}-\eta_hx'_{\bar h}L_j)=\operatorname{Tr}(\eta_h(L_ix_{\bar h}-x'_{\bar h}L_j))\]
Hence $L_ix_{\bar h}=x'_{\bar h}L_j$, and thus $\ima L$ is stable by $x'$. Moreover:\[
\forall i,\forall \tau_i, \operatorname{Tr}(\tau_if'_iL_i)=0\Rarrow\forall i, f'_iL_i=0\Rarrow\ima L\subset\ker f',\]
hence the lemma comes from the stability of $(x',f')$.
\edemo

We have to check that $V=V(X,F,G)$. It is easy to see that $V\subset V(X,F,G)$. Moreover:\[
F^{-1}(W)=\{v+v'\in V\oplus V'\mid f(v)+\theta(v')+f'(v')\in W\}=V\oplus\ker f',\]
hence, if $Y$ is an $X$-stable subspace of $F^{-1}(W)$, $Y/V$ is an $x'$-stable subspace of $\ker f'$. Since $(x',f')$ is stable, we have $Y\subset V$, and thus $V= V(X,F,G)$.

We have proved that the fiber $\textup{\sffamily T}_\circ^{-1}((x,f),(x',f'))$ is isomorphic to:\[
G_{\nu+\nu'}\times\CC^{\langle \lambda',\nu \rangle +(\nu',\nu)+\langle \nu',\lambda \rangle-\langle \nu',\nu\rangle}\]
and thus is connected. 
\edemo

\lemm\label{fibtensstab} Consider $(x,f,g)\in\widetilde\Z_\circ(\nu,\nu')$ and $ V= V(x,f,g)$. Then:\[
(x,f,g)\textup{ stable }\Lrarrow (x,f)_{| V\times  V}^{| V\times W}\textup{ stable}\]
and we denote by $\widetilde\Z(\nu,\nu')$ the subvariety of stable points of $\widetilde\Z_\circ(\nu,\nu')$, and:\[
\widetilde{\mathfrak Z}(\nu,\nu')=\widetilde\Z(\nu,\nu')/\!\!/G_{\nu+\nu'}.\]
\elemm

\demo 
The equivalence is a consequence of the following facts:\bite\itemsep0.5em
\item the restriction of a stable point is stable;
\item the extension of a stable point by a stable point is stable;
\item the point $(x,f)_{|(V_{\nu+\nu'}/ V)\times(V_{\nu+\nu'}/ V)}^{|(V_{\nu+\nu'}/ V)\times(W\oplus W'/W)}$ is stable.\eite\edemo

\theo\label{tenso}
We have the following bijection:
\[\xymatrix{
\Irr\mathfrak L(\nu,\lambda)\times\Irr\mathfrak L(\nu,\lambda')\ar[r]^-{\otimes}_-{\sim}&\Irr\widetilde{\mathfrak Z}(\nu,\nu').
}\]
\etheo

\demo
Define $\check\Z(\nu,\nu')$ as the variety of stable points of $\check\Z_\circ(\nu,\nu')$. We have the following diagram:\[
\xymatrix{
\check \Z(\nu,\nu')\ar[r]^-{\textup{\sffamily T}}\ar[d]& \La(\nu,\lambda)\times\La(\nu',\lambda')\ar[d]\\
\widetilde{\mathfrak Z}(\nu,\nu')\ar[r]^-{\mathfrak{T}}&\mathfrak L(\nu,\lambda)\times\mathfrak L(\nu',\lambda')  
}\]
where the rightmost vertical map is juste the free quotient by $G_\nu\times G_{\nu'}$. The leftmost map being a principal bundle with fibers isomorphic to:\[
G_\nu\times G_{\nu'}\times\textup{Grass}_{\nu,\nu'}^I(\nu+\nu')\times G_{\nu+\nu'},\]
we get our bijection thanks to~\ref{fibtens} and~\ref{fibtensstab}.
\edemo

Again, there is an alternative definition for $\widetilde{\mathfrak Z}(\nu,\nu')$, given in~\cite{nakatens}. 
Denote by $*$ the $\CC^*$-action on $\mathfrak M(\mathbf v,\lambda+\lambda')$ induced by the one parameter subgroup $\CC^*\rarrow GL(W\oplus W')$, $t\mapsto t\id_W\oplus\id_{W'}$. We have:\[
\mathfrak M(\mathbf v,\lambda+\lambda')^{\CC^*}\simeq\bigsqcup_{\nu+\nu'=\mathbf v}{\mathfrak M}(\nu,\lambda)\times\mathfrak M(\nu',\lambda')\]
and:\[
\widetilde {\mathfrak Z}(\nu,\nu')=\{[x,f,g]\in\mathfrak M(\mathbf v,\lambda+\lambda')\mid\lim_{t\to0}t*[x,f,g]\in\mathfrak L(\nu,\lambda)\times\mathfrak L(\nu',\lambda')\}.\]

Hence we also have, as in~\cite[3.15]{nakatens}, the following:
\prop The subvariety $\widetilde{\mathfrak Z}(\nu,\nu')\subset\mathfrak M(\nu+\nu',\lambda+\lambda')$ is Lagrangian. \eprop

The results of the section 2.2 lead to the following, completing~\cite[4.3]{nakatens} which deals with the case $\omega_i=0$:
\prop\label{bijtens}
Consider $i$ such that $\omega_i>0$ and $l>0$. If:\[
\lambda_i+\lambda'_i+\sum_{h\in H_i} \mathbf v_{t(h)}>0,\]
 we have a bijection:
\becen$\Irr\widetilde{\mathfrak Z}(\mathbf v)_{i,l}\isom\Irr\widetilde{\mathfrak Z}(\mathbf v-le_i)_{i,0}\times\Irr\Lambda({le_i})$.\ecen
\eprop

\subsection{Comparison of two crystal-type structures}

\nota For every $X\in\Irr\widetilde{\mathfrak Z}(\mathbf v)_{i,l}$, we will denote by $\eee_i(X)\in\Irr\Lambda(le_i)$ the composition of the second projection with the bijection obtained in~\ref{bijtens}, and $|\eee_i(X)|=l$. Note that if $(X,X')\in\Irr\mathfrak L(\nu,\lambda)\times\Irr\mathfrak L(\nu',\lambda')$, the quantity $\eee_i(X\otimes X')$ makes sense thanks to~\ref{tenso} and~\ref{bijtens}.\enota

We will write $\Omega(i)=\{b_{i,j}\}_{1≤j≤\omega_i}$ for $i$ imaginary, or $\Omega(i)=\{b_{j}\}_{1≤j≤\omega_i}$ if it is not ambiguous.

\lemm\label{lemmcyc} Let $i$ be an imaginary vertex and assume $\sum_{h\in H_i}n_{t(h)}>0$. For every $C\in\Irr\mathfrak L(\nu,\lambda)$, there exists $(x,f)\in C$, $v\in\ima\sum_{h\in H_i}x_{\bar h}$ such that:\[
\CC[ x_{\bar b_1}].v=\mathfrak I_i(x,f).\]
\elemm

\demo We proceed by induction on $\nu_i$, the first step being trivial. For the inductive step, we can immediatly conclude if $C\in\Irr\mathfrak L(\nu,\lambda)_{i,l}$ for $l>0$. Otherwise, $C\in\Irr\mathfrak L(\nu,\lambda)_{i,0}$, but $C\in\Irr\mathfrak L(\nu,\lambda)_{j,l}$ for some $j\in I$ and $l>0$. 
There exists a minimal chain $(j_k,l_k,C_k)_{1≤k≤s}$ of elements of $I\times\NN_{>0}\times\Irr\mathfrak L(-,\lambda)$ such that:\bite
\item $(j_1,l_1,C_1)=(j,l,C)$;
\item $C_{k+1}=\textup{pr}_1\mathfrak l(\nu-l_{1}j_{1}-\dots-l_kj_k,\lambda)_{j_k,l_k}(C_{k})$ where $\textup{pr}_1$ is the first projection;
\item $j_s=i$.\eite
We necessarely have $j_{s-1}$ adjacent to $i$, and by the induction hypothesis, the proposition is satisfied by $C_s$, and thus by $C_{s-1}$. But then, thanks to~\ref{lemm0} and~\ref{lemm1}, the proposition is also satisfied by $C_{s-2}$ for a generic choice of $\eta_{\bar h}$ (using the notations of the proof of~\ref{lemmclef} where $i$ is replaced by $j_{s-1}$). Hence it is also satisfied by $C=C_1$.
\edemo

\prop\label{critens} Let $i$ be an imaginary vertex and consider $(X, X')\in\Irr\mathfrak L(\nu,\lambda)\times\Irr\mathfrak L(\nu',\lambda')$. Assume $|\eee_i(X')|<\nu'_i$ or $0<\lambda'_i$. Then we have:\[
\eee_i(X\otimes X')=\eee_i(X').\]
\eprop

\demo Put $(Y,C)=\mathfrak l(n,m)_{i,l}(X)$ where $l=|\eee_i(X)|$. Take $((x,f),(x',f'))\in X\times X'$. Consider the equation $\zeta_i=0$ used in the proof of~\ref{lemmfibtens}:\[
\tau_if'_i+\dsum_{h\in H:s(h)=i}\epsilon(\bar h)(x_{\bar h}\eta_h+\eta_{\bar h}x'_h)=0.\]
Note $\eta_{b_j}=\eta_j$, $x_{b_j}=x_j$ and $x_{\bar b_j}=\bar x_j$ (and the same with $x'$), take $\eta_{\bar b_j}=0$ so that our equation becomes:\begin{align*}
\tau_if'_i+\dsum_{h\in H_i}\eta_{\bar h}x'_h&=\dsum_{1≤j≤\omega_i}(\bar x_j\eta_j-\eta_j {\bar x'_j})\\
&=\bar x_1\eta_1-\eta_1 {\bar x'_1}\end{align*}
if we also set $\eta_j=0$ for $j≥2$ (if any). Then, we set:\begin{align*}
x'&=f'_i+\doplus_{h\in H_i}x'_h:V_{\nu'_i}\rarrow W'_i\oplus\doplus_{h\in H_i}V_{\nu'_{t(h)}}\\
\bar\eta&=\tau_i+\dsum_{h\in H_i}\eee(\bar h)\eta_{\bar h}:W'_i\oplus\doplus_{h\in H_i}V_{\nu'_{t(h)}}\rarrow V_{\nu_i}\\
\bar x&=\dsum_{h\in H_i}\eee(\bar h)x_{\bar h}:\doplus_{h\in H_i}V_{\nu_{t(h)}}\rarrow V_{\nu_i}\\
\eta&=\doplus_{h\in H_i}\eta_{h}:V_{\nu'_i}\rarrow\doplus_{h\in H_i}V_{\nu_{t(h)}}
\end{align*}
and our equation finally becomes:\becen$
\bar\eta x'+\eta\bar x=\bar x_1\eta_1-\eta_1 {\bar x'_1}$.\ecen
Consider the open subvariety of $X\times X'$ where:\benu
\item there exists $\mathbf v\in V_{\nu_i}$ such that its image $\bar{\mathbf v}\in V_{\nu_i}/\mathfrak I_i(x,f)$ satisfies:\[
\CC[ \bar x_{1|V_{\nu_i}/\mathfrak I_i(x,f)}].\bar{\mathbf v}=V_{\nu_i}/\mathfrak I_i(x,f);\]
\item $\bar x'_1$, $\bar x_{1|\mathfrak I_i(x,f)}$ and $\bar x_{1|\CC^{n_i}/\mathfrak I_i(x,f)}$ have disjoint spectra;
\item there exist $v$ and $v'$ such that $\mathbf w=\sum_{h\in H_i}x_{\bar h}(v)$ and $\mathbf w'=\sum_{h\in H_i}x'_{\bar h}(v')$ satisfy:\[
\CC[ \bar x_1\oplus \bar x'_1].\mathbf w\oplus\mathbf w'=\mathfrak I_i(x,f)\oplus\mathfrak I_i(x',f');\]
\eenu
which is nonempty, thanks to~\ref{lemm1},~\ref{lemmcyc} and~\ref{lemm0}. Take:\bite
\item $\bar\eta=\tau_i$ and $\mathbf v\in\ima\tau_i$ if $\lambda'_i>0$;
\item $\bar \eta$ such that $\bar\eta(v')=\mathbf v$ if $\nu'_i>|\eee_i(X')|$ (possible since $v'\neq0$).\eite
From~\ref{lemm0}, we get (with the notations used in the proof of~\ref{fibtens}):\[
\CC[ X_{\bar b_1}].\ima\Big(\sum_{h\in H_i}X_{\bar h}\Big)=V_{\nu_i}\oplus\mathfrak I_i(x',f').\]
We have to check that we can choose $\eta$ such that the equations $\zeta_{t(h)}=0$ are satisfied for every $h\in H_i$ (if $\lambda'_i>0$ and $\bar\eta=\tau_i$, just take $\eta=0$). It suffices to set $\eta_hx'_{\bar h}(v'_{t(h)})=-x_h\eta_{\bar h}(v'_{t(h)})$ (possible since $\nu'_i>|\eee_i(X')|$ and since we may assume that $v'_{t(h)}=0$ if $x'_{\bar h}(v'_{t(h)})=0$) and to set $\eta$ and $\bar\eta$ equal to zero on supplementaries of $\CC\mathbf w'$ and $\CC v'$ respectively. 
We can finally choose $\eta_1$ such that $\bar\eta x'+\eta\bar x=\bar x_1\eta_1-\eta_1 {\bar x'_1}$ (possible since $\Spec\bar x'_1\cap\Spec\bar x_1=\varnothing$). 
Since:\[ 
\codim\mathfrak I_i(x,f)≥|\eee_i(X')|,\]
for every $(x,f)\in X\otimes X'$, the subvariety of $X\otimes X'$ defined by:\[
\codim\mathfrak I_i(x,f)=|\eee_i(X')|,\]
 is open, and we have shown it is non empty, hence the theorem is proved.
\edemo

\prop Assume $\lambda'_i=0$, $|\eee_i(X')|=\nu'_i$ and $\sum_{h\in H_i}\nu'_{t(h)}>0$. Then we still have $\eee_i(X\otimes X')=\eee_i(X')$.\eprop

\demo
Thanks to the previous proof, the result is clear if there exists an imaginary vertex $j$ adjacent to $i$: the choice of $x_{\bar b_{j,1}}$ and $x'_{\bar b_{j,1}}$ with disjoint spectra enables to use $\eta_{b_{j,1}}$ for $\zeta_{j}=0$ to be satisfied (with the usual notation $\Omega(j)=\{b_{j,1},\ldots,b_{j,\omega_j}\}$).

Assume that every neighbour of $i$ is real. Following the previous proof, assume $\bar \eta=\eta_{\bar h}$ is of rank $1$ for some $h:i\rarrow j$. We have to check that $\zeta_j=0$ can be satisfied. It is clear if $f'_j\neq 0$: just choose $\tau_j$ such that $\tau_jf'_j=-\eee( h)x_h\eta_{\bar h}$ and $\eta_p=0=\eta_{\bar p}$ if $p\in H_j\setminus\{\bar h\}$, so that $\zeta_j=0$ is satisfied.  Otherwise, there necessarily exists an edge $q:j\rarrow k\neq i$ such that $x'_q\neq 0$ (if not, $V'_{\nu'_i}\oplus V'_{\nu'_j}\subseteq\ker f'$ would be $x'$-stable, which is not possible for every vertex $j$ adjacent  to $i$ since $\sum_{h\in H_i}\nu'_{t(h)}>0$).
Hence it is possible to choose $\eta_{\bar q}$ so that $\eee(\bar q)\eta_{\bar q}x'_q=-\eee(h)x_h\eta_{\bar h}$ and $\eta_p=0=\eta_{\bar p}$ if $p\in H_j\setminus\{\bar h,q\}$, and thus get $\zeta_j=0$ satisfied.
\edemo

We have proved the following:

\theo\label{tenscrysgeom}
Let $i$ be an imaginary vertex and consider $(X, X')\in\Irr\mathfrak L(\nu,\lambda)\times\Irr\mathfrak L(\nu',\lambda')$. We have:\[
\eee_i(X\otimes X')=\left\{\begin{aligned}&\eee_i(X')&&\textup{if }\lambda'_i+\sum_{h\in H_i}\nu'_{t(h)}>0\\
&\eee_i(X)&&\textup{otherwise.}\end{aligned}\right.\]
\etheo

\section{Generalized crystals}

Let $(-,-)$  denote the symmetric Euler form on $\ZZ I$: $(i,j)$ is equal to the opposite of the number of edges of $\Omega$ between $i$ and $j$ for $i\neq j\in I$, and $(i,i)=2-2\omega_i$. We will still denote by $I^\text{re}$ (resp. $I^\text{im}$) the set of real (resp. imaginary) vertices, and by $I^\text{iso}\subseteq I^\text{im}$ the set of \emph{isotropic} vertices: vertices $i$ such that $(i,i)=0$, \textit{i.e.}\ such that $\omega_i=1$.
We also set $I_\infty=(I^\text{re}\times\{1\})\cup( I^\text{im}\times\NN_{≥1})$, and $(\iota,j)=l(i,j)$ if $\iota=(i,l)\in I_\infty$ and $j\in I$.

\subsection{A generalized quantum group}

We recall some of the definitions and results exposed in~\cite[§2]{article2}.

\defi Let {\sffamily F} denote the  $\QQ(v)$-algebra generated by $(E_\iota)_{\iota\in I_\infty}$, naturally $\NN I$-graded by $\text{deg}(E_{i,l})=li$ for $(i,l)\in I_\infty$. We put $\text {\sffamily F}[A]=\{x\in\textup{\sffamily F}\mid |x|\in A\}$ for any $A\subseteq\NN I$, where we denote by $|x|$ the degree of an element $x$.\edefi

For $\nu=\sum\nu_ii\in\ZZ I$, we set:\benu
\item[$\triangleright$] ht$(\nu)=\sum\nu_i$ its height;
\item[$\triangleright$] $v_\nu=\prod v_i^{\nu_i}$ if $v_i=v^{(i,i)/2}$.
\eenu
We endow $\text{\sffamily F}\otimes\textup{\sffamily F}$ with the following multiplication:\[
(a\otimes b)(c\otimes d)=v^{(|b|,|c|)}(ac)\otimes(bd).\]
and equip {\sffamily F} with a comultiplication $\delta$ defined by:\begin{align*}
\delta(E_{i,l})&=\dsum_{t+t'=l}v_i^{tt'}E_{i,t}\otimes E_{i,t'}\end{align*}
where $(i,l)\in I_\infty$ and $E_{i,0}=1$.

\prop For any family $(\nu_\iota)_{\iota\in I_\infty}$, we can endow \textup{\sffamily F} with a bilinear form $\{-,-\}$ such that:\benu
\item[$\triangleright$]$\{ x,y \}=0$ if $|x|\neq|y|$;
\item[$\triangleright$]$\{ E_\iota,E_\iota\}=\nu_\iota$ for all $\iota\in I_\infty$;
\item[$\triangleright$]$\{ ab,c\}=\{ a\otimes b,\delta (c)\}$ for all $a,b,c\in\textup{\sffamily F}$.
\eenu
\eprop

\nota Take $i\in I^\text{im}$ and $\textup{\sffamily c}$ a composition or a partition. We put $E_{i,\textup{\sffamily c}}=\prod_j E_{i,\textup{\sffamily c}_j}$ and $ \nu_{i,\textup{\sffamily c}}=\prod_j\nu_{i,\textup{\sffamily c}_j}$.  If $i$ is real, we will often use the index $i$ instead of $i,1$.
\enota

\prop Consider $(\iota,j)\in I_\infty\times I^\text{re}$. The element:\begin{align}\label{Serre}
\dsum_{t+t'=-(\iota,j)+1}(-1)^tE_j^{(t)}E_{\iota}E_j^{(t')}\end{align}
belongs to the radical of $\{-,-\}$.\eprop

\defi We denote by $\tilde U^+$ the quotient of \textup{\sffamily F} by the ideal spanned by the elements~\ref{Serre} and the commutators $[E_{i,l},E_{i,k}]$ for every isotropic vertex $i$, so that $\{-,-\}$ is still defined on $\tilde U^+$. We denote by $ U^+$ the quotient of $\tilde U^+$ by the radical of $\{-,-\}$.
 \edefi

\defi
Let $\hat U$ be the quotient of the algebra generated by $K_i^\pm$, $E_\iota$, $F_\iota$ ($i\in I$ and $\iota\in I_\infty$) subject to the following relations:\[
K_iK_j&=K_jK_i\\
K_iK_i^-&=1\\
K_jE_\iota&=v^{(j,\iota)}E_\iota K_j\\
K_jF_\iota&=v^{-(j,\iota)}F_\iota K_j\\
\sum_{t+t'=-(\iota,j)+1}(-1)^tE_j^{(t)}E_{\iota}E_j^{(t')}&=0\qquad(j\in I^\text{re})\\
\sum_{t+t'=-(\iota,j)+1}(-1)^tF_j^{(t)}F_{\iota}F_j^{(t')}&=0\qquad(j\in I^\text{re})\\
[E_{i,l},E_{j,k}]&=0\qquad\text{if }(i,j)=0\\
[F_{i,l},F_{j,k}]&=0\qquad\text{if }(i,j)=0\\
[E_{i,l},E_{i,k}]&=0\qquad(i\in I^\text{iso})\\
[F_{i,l},F_{i,k}]&=0\qquad(i\in I^\text{iso}).\]
We extend the graduation by $|K_i|=0$ and $|F_\iota|=-|E_\iota|$, and we set $K_\nu=\prod_iK_i^{\nu_i}$ for every $\nu\in\ZZ I$.

We endow $\hat U$ with a comultiplication $\Delta$ defined by:\[
\Delta(K_i)&=K_i\otimes K_i\\
\Delta(E_{i,l})&=\dsum_{t+t'=l}v_i^{tt'}E_{i,t}K_{t'i}\otimes E_{i,t'}\\
\Delta(F_{i,l})&=\dsum_{t+t'=l}v_i^{-tt'}F_{i,t}\otimes K_{-ti}F_{i,t'}.\]

We extend $\{-,-\}$ to the subalgebra $\hat U^{≥0}\subseteq\hat U$ spanned by $(K_i^\pm)_{i\in I}$ and $(E_\iota)_{\iota\in I_\infty}$ by setting $\{ xK_i,yK_j\}=\{ x,y\} v^{(i,j)}$ for $x,y\in\tilde U^+$.

We use the Drinfeld double process to define $\tilde U$ as the quotient of $\hat U$ by the relations:\begin{align}\label{DD}
\dsum \{ a_{(1)},b_{(2)}\} \omega(b_{(1)})a_{(2)}&=\dsum \{ a_{(2)},b_{(1)}\} a_{(1)}\omega(b_{(2)})
\end{align}
for any $a,b\in \tilde U^{≥0}$, where $\omega$ is the unique involutive automorphism of $\hat U$ mapping $E_\iota$ to $F_\iota$ and $K_i$ to $K_{-i}$, and where we use the Sweedler notation, for example $\Delta(a)=\sum a_{(1)}\otimes a_{(2)}$.

Setting $x^-=\omega(x)$ for $x\in \tilde U$, we define $\{-,-\}$ on the subalgebra $\tilde U^{-}\subseteq\tilde U$ spanned by $(F_\iota)_{\iota\in I_\infty}$ by setting $\{ x,y\}=\{ x^-,y^-\}$ for any $x,y\in\tilde  U^-$. We will denote by $U^-$ (resp. $U$) the quotient of $\tilde U^-$ (resp. $\tilde U$) by the radical of $\{-,-\}$ restricted to $\tilde U^-$(resp. restricted to $\tilde U^-\times\tilde U^+$).\edefi

\prop\label{hypo} Assume:\begin{align*}
\{ E_{\iota},E_{\iota}\}\in1+v^{-1}\NN[\![ v^{-1}]\!].\end{align*}
for every $\iota\in I_\infty$. Then we have $\tilde U^-\simeq U^-$.\eprop

\nota We denote by $\mathcal C_{i,l}$ the set of compositions $\textup{\sffamily c}$ (resp. partitions) such that $|\textup{\sffamily c}|=l$ if $(i,i)<0$ (resp. $(i,i)=0$), and $\mathcal C_i=\sqcup_{l≥0}\mathcal C_{i,l}$. If $i$ is real, we just put $\mathcal C_{i,l}=\{l\}$.

Denote by $u\mapsto\bar u$ the involutive $\QQ$-morphism of $U$ stabilizing $E_\iota$, $F_\iota$, and mappinf $K_i$ to $K_{-i}$, and $v$ to $v^{-1}$.
\enota

\prop\label{prim} For any imaginary vertex $i$ and any $l≥1$, there exists a unique element $a_{i,l}\in U^+[li]$ such that, if we set $b_{i,l}=a_{i,l}^-$, we get:\benu
\item $\QQ(v)\langle E_{i,l}\mid l≥1\rangle=\QQ(v)\langle a_{i,l}\mid l≥1\rangle$ and $\QQ(v)\langle F_{i,l}\mid l≥1\rangle=\QQ(v)\langle b_{i,l}\mid l≥1\rangle$ as algebras;
\item $\{ a_{i,l},z\}=\{ b_{i,l},z^-\}=0$ for any $z\in\QQ(v)\langle E_{i,k}\mid k<l\rangle$;
\item $a_{i,l}-E_{i,l}\in\QQ(v)\langle E_{i,k}\mid k<l\rangle$ and $b_{i,l}-F_{i,l}\in\QQ(v)\langle F_{i,k}\mid k<l\rangle$;
\item $\bar a_{i,l}=a_{i,l}$ and $\bar b_{i,l}=b_{i,l}$;
\item $\delta(a_{i,l})=a_{i,l}\otimes1+1\otimes a_{i,l}$ and $\delta(b_{i,l})=b_{i,l}\otimes 1+1\otimes b_{i,l}$;
\eenu
\eprop

\nota Consider $i\in I^\text{im}$ and $\textup{\sffamily c}\in\mathcal C_{i,l}$. We set $\tau_{i,l}=\{ a_{i,l},a_{i,l}\}$,  $a_{i,\textup{\sffamily c}}=\prod_j a_{i,\textup{\sffamily c}_j}$, and $ \tau_{i,\textup{\sffamily c}}=\prod_j\tau_{i,\textup{\sffamily c}_j}$. Notice that $\{a_{i,\textup{\sffamily c}}\mid \textup{\sffamily c}\in\mathcal C_{i,l}\}$ is a basis of $U^+[l i]$.
\enota

\defi
We denote by $\delta_{i,\textup{\sffamily c}},\delta^{i,\textup{\sffamily c}}:U^+\rarrow U^+$ the linear maps defined by: \[
\delta(x)&=\dsum_{\textup{\sffamily c}\in\mathcal C_{i,l}}\delta_{i,\textup{\sffamily c}}(x)\otimes a_{i,\textup{\sffamily c}}+\text{obd}\\
\delta(x)&=\dsum_{\textup{\sffamily c}\in\mathcal C_{i,l}}a_{i,\textup{\sffamily c}}\otimes\delta^{i,\textup{\sffamily c}}(x)+\text{obd} \]
where "obd" stands for terms of bidegree not in $\NN I\times\NN i$ in the former equality, $\NN i\times\NN I$ in the latter one.
\edefi

\subsection{Kashiwara operators}

\prop\label{kashiprop}
Let $i$ be an imaginary vertex, $l>0$, $\textup{\sffamily c}=(\textup{\sffamily c}_1,\ldots,\textup{\sffamily c}_r)\in\mathcal C_{i}$ and $(y,z)\in (U^+)^2$. We have the following identities:\benu\itemsep0.5em
\item $\delta^{i,l}(yz)=\delta^{i,l}(y)z+v^{l(i,|y|)}y\delta^{i,l}(z)$;
\item $[a_{i,l},z^-]=\tau_{i,l}\big\{\delta_{i,l}(z^-)^-K_{-li}-K_{li}\delta^{i,l}(z^-)^-\big\}$;
\item $\delta^{i,l}(a_{i,\textup{\sffamily c}})=\dsum_{k:\textup{\sffamily c}_k=l}v_i^{2l\textup{\sffamily c}_{k-1}}a_{i,\textup{\sffamily c}\backslash\textup{\sffamily c}_k}$
\eenu
where $\textup{\sffamily c}_{0}=0$ and $\textup{\sffamily c}\backslash\textup{\sffamily c}_k=(\textup{\sffamily c}_1,\ldots,\hat{\textup{\sffamily c}}_k,\ldots,\textup{\sffamily c}_r)$, the notation $\hat{\textup{\sffamily c}}_k$ meaning that $\textup{\sffamily c}_k$ is removed from $\textup{\sffamily c}$. 
\eprop

\demo The first equality comes from the definition of $\delta^{i,l}$, the second from the primitive character of $a_{i,l}$ and the formula~\ref{DD} with $a=a_{i,l}$ and $b=z^-$. The third comes from the definition of $\delta_{i,l}$ and the primitive character of the $a_{i,h}$.\edemo

\defi Define $e'_{i,l}:U^-\rarrow U^-$ by $e'_{i,l}(z^-)=\delta^{i,l}(z)^-$ for any $z\in U^+$.\edefi

\prop Set \[\mathcal K_i=\dcap_{l>0}\ker e'_{i,l}\] for any $i\in I^\textup{im}$. We have the following decomposition:\[
U^-=\bigoplus_{\textup{\sffamily c}\in\mathcal C_i}b_{i,\textup{\sffamily c}}\mathcal K_i.\]
\eprop

\demo Let's first prove the existence. Consider $u\in U^-$, and assume first that $u$ is of the following form: $u=mb_{i,\textup{\sffamily c}}m'$ for some $\textup{\sffamily c}\in\mathcal C_i$ and some $m,m'\in\mathcal K_i$. We proceed by induction on $|\textup{\sffamily c}|$. If $|\textup{\sffamily c}|=0$, we have $mm'\in\mathcal K_i$ thanks to~\ref{kashiprop} (1). Otherwise, set $[y,z]^\circ=v^{-(|y|,|z|)}yz-zy$ for any $y,z\in U^+$. Thanks to~\ref{kashiprop} (1), and since $\delta^{i,l}(a_{i,k})=\delta_{l,k}$, we have for any $y\in\cap_{l>0}\ker\delta^{i,l}$ and any $k>0$:\[
\delta^{i,l}([y,a_{i,k}]^\circ)&=v^{-k(i,|y|)}\delta^{i,l}(ya_{i,k})-\delta^{i,l}(a_{i,k}y)\\
&=v^{-k(i,|y|)}v^{l(i,|y|)}y\delta^{i,l}(a_{i,k})-\delta^{i,l}(a_{i,k})y\\
&=\delta_{l,k}v^{(l-k)(i,|y|)}y-\delta_{l,k}y\\
&=0.\]
Hence, the following equality:\[
u=v^{\textup{\sffamily c}_1(|m|,i)}[m,b_{i,\textup{\sffamily c}_1}]^\circ b_{i,\textup{\sffamily c}\backslash\textup{\sffamily c}_1}m'+v^{-\textup{\sffamily c}_1(|m|,i)}b_{i,\textup{\sffamily c}_1}mb_{i,\textup{\sffamily c}\backslash\textup{\sffamily c}_1}m'
\]
along with the induction hypothesis allow us to conclude since $|\textup{\sffamily c}\backslash\textup{\sffamily c}_1|<|\textup{\sffamily c}|$, and since $\oplus_{\textup{\sffamily c}\in\mathcal C_i}b_{i,\textup{\sffamily c}}\mathcal K_i$ is stable by left-multiplication by $b_{i,\textup{\sffamily c}_1}$.

Then, we prove the existence of the decompostion for a general $u\in U^-$, using induction on $-|u|$. If $u\neq1$, we can write:\[
u=\dsum_{\iota\in I_\infty}b_{\iota}u_\iota\]
for some finitely many nonzero $u_\iota\in U^-$. Thanks to our induction hypothesis, we have:\[
u=\dsum_{\iota\in I_\infty,\textup{\sffamily c}\in\mathcal C_i}b_\iota b_{i,\textup{\sffamily c}} z_{\iota,\textup{\sffamily c}}\]
for some finitely many nonzero $z_{\iota,\textup{\sffamily c}}\in\mathcal K_i$. Then:\[
u=\dsum_{l>0,\textup{\sffamily c}\in\mathcal C_i}b_{i,(l,\textup{\sffamily c})}z_{(i,l),\textup{\sffamily c}}+\dsum_{\substack{\iota\in I_\infty\backslash(\{i\}\times\NN_{>0})\\\textup{\sffamily c}\in\mathcal C_i}}b_\iota b_{i,\textup{\sffamily c}} z_{\iota,\textup{\sffamily c}}\]
and we have the result since $b_\iota b_{i,\textup{\sffamily c}} z_{\iota,\textup{\sffamily c}}$ is of the form $mb_{i,\textup{\sffamily c}}m'$ for some $m,m'\in\mathcal K_i$. Indeed, it is straightforward from the definitions that $\delta^{i,l}(a_{j,h})=0$ for any $l,h>0$ if $j\neq i$. Note that if $i\notin I^\textup{iso}$, the composition $(l,\textup{\sffamily c})$ is the composition $\textup{\sffamily c}'$ where $\textup{\sffamily c}'_1=l$ and $\textup{\sffamily c}'_k=\textup{\sffamily c}_{k-1}$ if $k≥2$, but if $i\in I^\textup{iso}$, $(l,\textup{\sffamily c})$ stands for the partition $\textup{\sffamily c}\cup l$.

To prove the unicity of the decomposition, consider a minimal nontrivial relation of dependance:\[
0=\dsum_{\textup{\sffamily c}\in\mathcal C_i}a_{i,\textup{\sffamily c}}z_\textup{\sffamily c},\]
where $z_\textup{\sffamily c}\in\mathcal K_i^-$. We have to considerate separately the cases $i\in I^\textup{iso}$ and $i\notin I^\textup{iso}$.
First, consider $i\notin I^\textup{iso}$.
Consider $r$ maximal such that there exists $\textup{\sffamily c}=(\textup{\sffamily c}_1,\ldots,\textup{\sffamily c}_r)$ such that $z_\textup{\sffamily c}\neq0$.
Using~\ref{kashiprop} (1) and applying repeatedly~\ref{kashiprop} (3), we see that for any $\textup{\sffamily c}'\in\mathfrak S_r\textup{\sffamily c}$ (with the convention $(\sigma\textup{\sffamily c})_k=\textup{\sffamily c}_{\sigma(k)}$):\[
0&=\delta^{i,\textup{\sffamily c}'}\bigg(\dsum_{\textup{\sffamily c}''\in\mathcal C_i}a_{i,\textup{\sffamily c}''}z_{\textup{\sffamily c}''}\bigg)\\
&=\dsum_{\textup{\sffamily c}''\in\mathcal C_i}\delta^{i,\textup{\sffamily c}'}(a_{i,\textup{\sffamily c}''})z_{\textup{\sffamily c}''}\\
&=\dsum_{\textup{\sffamily c}''\in\mathfrak S_r\textup{\sffamily c}}\delta^{i,\textup{\sffamily c}'}(a_{i,\textup{\sffamily c}''})z_{\textup{\sffamily c}''}\\
&=\dsum_{\textup{\sffamily c}''\in\mathfrak S_r\textup{\sffamily c}}P_{\textup{\sffamily c}',\textup{\sffamily c}''}(v_i)z_{\textup{\sffamily c}''}\]
where $P_{\textup{\sffamily c}',\textup{\sffamily c}''}(v)\in\ZZ[v]$. The third is equality is true by maximality of $r$. Since $(z_{\textup{\sffamily c}''})_{\textup{\sffamily c}''\in\mathfrak S_r\textup{\sffamily c}}\neq0$, we have to prove that:\[
\Delta(v)=\det(P_{\textup{\sffamily c}',\textup{\sffamily c}''}(v))_{\textup{\sffamily c}',\textup{\sffamily c}''\in\mathfrak S_r\textup{\sffamily c}}\neq 0\in\ZZ[v]\]
to end our proof in the case $(i,i)<0$ (since then we have $v_i\neq1$).
But, for any $\textup{\sffamily c}'\in\mathfrak S_r\textup{\sffamily c}$, one has, using~\ref{kashiprop} (3):\[
\underset{\textup{\sffamily c}''\in\mathfrak S_r\textup{\sffamily c}}{\textup{max}}\{\textup{deg}(P_{\textup{\sffamily c}',\textup{\sffamily c}''})\}=\textup{deg}(P_{\textup{\sffamily c}',{\textup{\sffamily c}}'})=\dsum_{1≤k≤r}\textup{\sffamily c}_{k-1}\textup{\sffamily c}_k=m\]
which is only reached for $\textup{\sffamily c}''={\textup{\sffamily c}}'$. Note that if $m=0$, the initial relation of dependance can be written:\[
0=\dsum_{l>0}a_{i,l}z_l,\]
 and we get $z_l=0$ after applying $\delta^{i,l}$. Otherwise, note that the degree of $\Delta$ is $|\mathfrak S_r\textup{\sffamily c}|m>0$, and in particular $\Delta\neq0$.

We finally have to prove the uniqueness in the case $(i,i)=0$. Write a relation of dependance of minimal degree:\[
0=\dsum_{\lambda\in\mathcal C_i}a_{i,\lambda}z_\lambda,\]
where $z_\lambda\in\mathcal K_i^-$. For any $\lambda$ and $l>0$, set $m_l(\lambda)=|\{s:\lambda_s=l\}|$, and denote by $\lambda\backslash l$ the partition obtained removing one of the $\lambda_s=l$ when $m_l(\lambda)≥1$. Hence $m_l(\lambda\backslash l)=m_l(\lambda)-1$. We have, thanks to~\ref{kashiprop} (1,3):\[
\delta^{i,l}\bigg(\dsum_{\lambda\in\mathcal C_i}a_{i,\lambda}z_\lambda\bigg)&=\dsum_{\lambda\in\mathcal C_i}m_l(\lambda)a_{i,\lambda\backslash l}z_\lambda\\
&=\dsum_{\mu\in\mathcal C_i}a_{i,\mu}\big\{(m_l(\mu)+1)z_{\mu\cup l}\big\}\]
which contradicts the minimality of the first relation. Note that the proof is easier in this case because we are dealing with partitions, hence the quantity $\mu\cup l$ is "uniquely defined".
\edemo

The following definition generalizes the Kashiwara operators (see \textit{e.g.}~\cite[2.3.1]{kashisaito}):
\defi\label{kashiop}
If $i$ is imaginary and $z=\sum_{\textup{\sffamily c}\in\mathcal C_i}b_{i,\textup{\sffamily c}}z_\textup{\sffamily c}\in U^-$, set:\[
\tilde e_{i,l}(z)&=\left\{\begin{aligned}&\dsum_{\textup{\sffamily c}:\textup{\sffamily c}_1=l}b_{i,\textup{\sffamily c}\backslash\textup{\sffamily c}_1}z_\textup{\sffamily c}&&\text{if }i\notin I^\text{iso}\\
&\dsum_{\lambda\in\mathcal C_i}\sqrt{\frac{m_l(\lambda)}{l}}b_{i,\lambda\backslash l}z_\lambda&&\text{if }i\in I^\text{iso}\end{aligned}\right.\\
\tilde f_{i,l}(z)&=\left\{\begin{aligned}&\dsum_{\textup{\sffamily c}\in\mathcal C_i}b_{i,(l,\textup{\sffamily c})}z_\textup{\sffamily c}&&\text{if }i\notin I^\text{iso}\\
&\dsum_{\lambda\in\mathcal C_i}\sqrt{\frac{l}{m_l(\lambda)+1}}b_{i,\lambda\cup l}z_\lambda&&\text{if }i\in I^\text{iso}.\end{aligned}\right.\]
\edefi

\subsection{Definition of generalized crystals}

Denote by $P$ the lattice $\ZZ^{I}$, still endowed with the pairing $\langle-,-\rangle$ defined by $\langle e_i,e_j\rangle=\delta_{i,j}$, where $e_i=(\delta_{i,j})_{j\in I}$ for every $i\in I$. We will also note $\aaa_i$ instead of $Ce_i$ where $C=((i,j))_{i,j\in I}$ still denotes the Cartan matrix associated to $Q$.

\defi\label{axiocrys} We call $Q$-crystal a set $\mathcal B$ together with maps:\[
 \textup{wt}&:\mathcal B\rarrow P\\
 \eee_i&:\mathcal B\rarrow\mathcal C_i\sqcup\{-\infty\}\\
 \phi_i&:\mathcal B\rarrow\NN\sqcup\{+\infty\}\\
 \tilde{ e}_{i},\tilde{ f}_{i}&:\mathcal B\rarrow\mathcal B\sqcup\{0\}&&i\in I^\textup{re}\\
 \tilde{ e}_{i,l},\tilde{ f}_{i,l}&:\mathcal B\rarrow\mathcal B\sqcup\{0\}&&i\in I^\textup{im}, l>0\]
 such that for every $b,b'\in\mathcal B$:\benu\itemsep0.5em
 \item[(A1)] $\langle e_i,\text{wt}(b)\rangle≥0$ if $i\in I^\text{im}$;
\item[(A2)]	$\textup{wt}(\tilde{ e}_{i,l}b)=\textup{wt}(b)+l\aaa_i$ if $\tilde{ e}_{i,l}b\neq0$;
\item[(A3)] $\textup{wt}(\tilde{ f}_{i,l}b)=\textup{wt}(b)-l\aaa_i$ if $\tilde{ f}_{i,l}b\neq0$;
\item[(A4)] $\tilde{ f}_{i,l}b=b'\Lrarrow b=\tilde{ e}_{i,l}b'$;
\item[(A5)] if $\tilde{ e}_{i,l}b\neq0$, $\eee_i(\tilde{ e}_{i,l}b)=\left\{\begin{aligned}&\eee_i(b)-l&&\textup{if }i\in I^\textup{re}\\
&\eee_i(b)\setminus l&&\textup{if }i\in I^\textup{im}\setminus I^\textup{iso}\textup{ and }l=\eee_i(b)_1\\
&\eee_i(b)\setminus l&&\textup{if }i\in I^\textup{iso};\end{aligned}\right.$
\item[(A6)] if $\tilde{ f}_{i,l}b\neq0$, $\eee_i(\tilde{ f}_{i,l}b)=\left\{\begin{aligned}&\eee_i(b)+l&&\textup{if }i\in I^\textup{re}\\
&(l,\eee_i(b))&&\textup{if }i\in I^\textup{im};\end{aligned}\right.$
\item[(A7)] $\phi_i(b)=\left\{\begin{aligned}&\eee_i(b)+\langle e_i,\textup{wt}(b)\rangle&&\textup{if }i\in I^\textup{re}\\
&+\infty&&\textup{if }i\in I^\textup{im}\textup{ and }\langle e_i,\textup{wt}(b)\rangle>0\\
&0&&\textup{otherwise,}\end{aligned}\right.$
\eenu
where, for $i\in I^\textup{re}$, we write $\tilde{ e}_{i,1},\tilde{ f}_{i,1}$ instead of $\tilde{ e}_i,\tilde{ f}_i$ and $\tilde{ e}_{i,l},\tilde{ f}_{i,l}$ instead of $\tilde{ e}^l_{i,1},\tilde{ f}^l_{i,1}$. Also, as earlier, $(l,\eee_i(b))$ stands for the partition $\eee_i(b)\cup l$ if $i\in I^\text{iso}$.

\edefi

\rema~

\bite
\item We will use the following notation: $\text{wt}_i=\langle e_i,\text{wt}\rangle$.
\item Note that this definition of $\phi_i$ already appears in~\cite{JKK}. Also note that since we will only be interested in normal crystals (see~\ref{definormal}), we require $|\eee_i|$ and $\phi$ to be non-negative, except if $\eee_i=-\infty$ in which case we set $|\eee_i|=-\infty$.
\item Set $\tilde{ e}_{i,\textup{\sffamily c}}=\tilde{ e}_{i,\textup{\sffamily c}_1}\ldots\tilde{ e}_{i,\textup{\sffamily c}_r}$ and $\tilde{ f}_{i,\textup{\sffamily c}}=\tilde{ f}_{i,\textup{\sffamily c}_1}\ldots\tilde{ f}_{i,\textup{\sffamily c}_r}$ for every $\textup{\sffamily c}=(\textup{\sffamily c}_1,\ldots,\textup{\sffamily c}_r)$. Set $\bar{\textup{\sffamily c}}=(\textup{\sffamily c}_r\ldots,\textup{\sffamily c}_1)$ if $\omega_i≥2$, $\bar{\textup{\sffamily c}}={\textup{\sffamily c}}$ if $\omega_i≤1$. We have:\[
\tilde{ f}_{i,\textup{\sffamily c}}b=b'\Lrarrow b=\tilde{ e}_{i,\bar{\textup{\sffamily c}}}b'.\]
\eite\erema

\exem For every vertex $i$, we define a crystal $\mathcal B_i$ by endowing $\mathcal C_i$ with the following maps:\[
\textup{wt}(\textup{\sffamily c})&=-|\textup{\sffamily c}|\aaa_i\\
\eee_i(\textup{\sffamily c})&=\textup{\sffamily c}\\
\eee_j(\textup{\sffamily c})&=-\infty\text{ if }j\neq i.\]
Note that the definition of $\tilde e_{i,l},\tilde f_{i,l}$ is dictated by the one of $\eee_i$, together with~\ref{axiocrys} (A5, A6). We will denote by $(0)_i$ the trivial element of $\mathcal C_i$.\eexem

\defi A morphism of crystals $\mathcal B_1\rarrow\mathcal B_2$ is a map $\mathcal B_1\sqcup\{0\}\rarrow\mathcal B_2\sqcup\{0\}$ mapping $0$ to $0$, preserving the weight, $\eee_i$, and commuting with the respective actions of the $\tilde e_\iota,\tilde f_{\iota}$ on $\mathcal B_1$ and $\mathcal B_2$.
\edefi

\defi\label{definormal} A crystal $\mathcal B$ is said to be \emph{normal} if for every $b\in\mathcal B$ and $i\in I$, we have:\[
\eee_i(b)&=\overline{\max\{\bar{\textup{\sffamily c}}\mid \tilde e_{i,{\textup{\sffamily c}}}(b)\neq0\}}\\
\phi_i(b)&=\max\{{|\textup{\sffamily c}}|\mid \tilde f_{i,{\textup{\sffamily c}}}(b)\neq0\}.\]
\edefi

\defi\label{defitensorprod} The \emph{tensor product} $\mathcal B\otimes\mathcal B'=\{b\otimes b'\mid b\in\mathcal B,b'\in\mathcal B'\}$ of two crystals is defined by:\benu\itemsep0.5em
\item $\textup{wt}(b\otimes b')=\textup{wt}(b)+\textup{wt}(b')$;
\item if $i\in I^\textup{re}$, $\eee_i(b\otimes b')=\max\big\{\eee_i(b'),\eee_i(b)-\textup{wt}_i(b')\big\}$;
\item if $i\in I^\textup{im}$, $\eee_i(b\otimes b')=\left\{\begin{aligned}&\eee_i(b')&&\textup{if }\phi_i(b)≥|\eee_i(b')|\\
&\eee_i(b)&&\textup{if }\phi_i(b)<|\eee_i(b')|;\end{aligned}\right.$
\item if $i\in I^\textup{re}$, $\phi_i(b\otimes b')=\max\big\{\phi_i(b')+\textup{wt}_i(b),\phi_i(b)\big\}$;
\item if $i\in I^\textup{im}$, $\phi_i(b\otimes b')=\left\{\begin{aligned}&\phi_i(b')&&\textup{if }\phi_i(b)≥|\eee_i(b')|\\
&\phi_i(b)&&\textup{if }\phi_i(b)<|\eee_i(b')|;\end{aligned}\right.$
\item for every $\iota=(i,l)\in I_\infty$, $\tilde{ e}_{\iota}(b\otimes b')=\left\{\begin{aligned} & b\otimes \tilde{ e}_{\iota}(b')&&\textup{if }\phi_i(b')≥|\eee_i(b)|\\
& \tilde{ e}_{\iota}(b)\otimes b'&&\textup{if }\phi_i(b')<|\eee_i(b)|;\end{aligned}\right.$
\item for every $\iota=(i,l)\in I_\infty$, $\tilde{ f}_{\iota}(b\otimes b')=\left\{\begin{aligned} & b\otimes \tilde{ f}_{\iota}(b')&&\textup{if }\phi_i(b')>|\eee_i(b)|\\
&\tilde{ f}_{\iota}(b)\otimes b'&&\textup{if }\phi_i(b')≤|\eee_i(b)|.\end{aligned}\right.$
\eenu
\edefi

\rema Note that when $i$ is imaginary, the condition $\phi_i(b')>|\eee_i(b)|$ is equivalent to:\becen
 $\phi_i(b')=+\infty$ or [$\phi_i(b')=0$ and $\eee_i(b)=-\infty$],\ecen
and $\phi_i(b')=|\eee_i(b)|$ is equivalent to $\phi_i(b')=0=|\eee_i(b)|$.
 \erema

\prop $\mathcal B\otimes\mathcal B'$ is a crystal if $\mathcal B'$ is normal.\eprop

\demo 
Note that the result is already known if $I^\text{im}=\varnothing$, hence we just have to check the axioms of~\ref{axiocrys} that concern imaginary vertices. Axioms (A1), (A2), (A3) and (A7) are clearly satisfied.

To prove that (A4) is satisfied, we first consider $b$ and $b'$ such that $\tilde e_{i,l}(b\otimes b')\neq0$. If $\phi_i(b')≥|\eee_i(b)|$, we have $\tilde e_{i,l}(b\otimes b')=b\otimes \tilde e_{i,l}(b')$. The crystal $\mathcal B'$ being normal, $\phi_i(\tilde e_{i,l}(b'))>0$ since $\tilde f_{i,l}\tilde e_{i,l}(b')=b'\neq0$. But $i$ is imaginary, so by definition $\phi_i(\tilde e_{i,l}(b'))\in\{0,+\infty\}$, and we get $\phi_i(\tilde e_{i,l}(b'))=+\infty$. Also by definition, $|\eee_i(b)|<+\infty$, hence we get:\[
\tilde f_{i,l}(b\otimes \tilde e_{i,l}b')=b\otimes \tilde f_{i,l}\tilde e_{i,l}b'=b\otimes b'.\]
If $\phi_i(b')<|\eee_i(b)|$, we have $\tilde e_{i,l}(b\otimes b')=\tilde e_{i,l}(b)\otimes b'$, where $|\eee_i(\tilde e_{i,l}(b))|=|\eee_i(b)|-l$ is necessarily nonnegative by definition of $\eee_i$. Also, $\phi_i(b')$ can not be equal to $+\infty$, hence is $0$, and:\[
\tilde f_{i,l}(\tilde e_{i,l}(b)\otimes b')=\tilde f_{i,l}(\tilde e_{i,l}(b))\otimes b'=b\otimes b'.\]
Assume now that  $\tilde f_{i,l}(b\otimes b')\neq0$. If $\phi_i(b')>|\eee_i(b)|$, we get $\tilde f_{i,l}(b\otimes b')=b\otimes \tilde f_{i,l}b'$. If $\eee_i(b)=-\infty$, we have $\phi_i(\tilde f_{i,l}b')>\eee_i(b)$. Otherwise, we necessarily have $\phi_i(b')=+\infty$. But:\[
\text{wt}_i(\tilde f_{i,l}(b'))=\text{wt}_i(b')-l\langle e_i,\aaa_i\rangle≥\text{wt}_i(b')\]
since $\langle e_i,\aaa_i\rangle≤0$ for every $i\in I^\text{im}$. Hence $\phi_i(\tilde f_{i,l}(b'))=+\infty$, and:\[
\tilde e_{i,l}(b\otimes \tilde f_{i,l}b')=b\otimes \tilde e_{i,l}\tilde f_{i,l}b'=b\otimes b'.\]
If $\phi_i(b')≤|\eee_i(b)|$, then $\tilde f_{i,l}(b\otimes b')=\tilde f_{i,l}(b)\otimes b'$, where:\[
|\eee_i(\tilde f_{i,l}(b))|=|\eee_i(b)|+l>|\eee_i(b)|≥\phi_i(b'),\]
hence:\[
\tilde e_{i,l}(\tilde f_{i,l}(b)\otimes b')=\tilde e_{i,l}(\tilde f_{i,l}(b))\otimes b'=b\otimes b'.\]

From the definitions and the proof of (A4) above, it is easy to check that:\[
\eee_i(\tilde e_{i,l}(b\otimes b'))=\left\{\begin{aligned} & \eee_i(\tilde e_{i,l}b')&&\text{if }\eee_i(b\otimes b')=\eee_i(b')\\
& \eee_i(\tilde e_{i,l}b)&&\text{if }\eee_i(b\otimes b')=\eee_i(b),\end{aligned}\right.
\]
hence (A5) is satisfied by $\mathcal B\otimes\mathcal B'$ since it is by $\mathcal B$ and $\mathcal B'$. For the same reasons, (A6) is satisfied, except if $\phi_i(b')=|\eee_i(b)|$, which can only happen if both are equal to $0$ (we still consider $i\in I^\text{im}$). But then $\tilde e_{i,l}(b')=0$, so there is nothing to prove. Otherwise we would have $\tilde f_{i,l}\tilde e_{i,l}b'=b'\neq0$, hence $\phi_i(\tilde e_{i,l}b')=+\infty$ by normality. But then:\[
\text{wt}_i(b')=\text{wt}_i(\tilde e_{i,l}(b'))-l\langle e_i,\aaa_i\rangle≥\text{wt}_i(\tilde e_{i,l}(b'))>0\]
would imply $\phi_i(b')=+\infty$ which contradicts the assumption.
\edemo

\subsection{The crystal $\mathcal B(\infty)$}

\subsubsection{Algebraic definition}
Let $\mathcal A\subset\QQ(v^{-1})$ be the subring consisting of rational functions without pole at $v^{-1}=0$, and $\mathcal L(\infty)$ be the sub-$\mathcal A$-module of $U^-$ generated by the elements $\tilde f_{\iota_1}\ldots\tilde f_{\iota_s}.1$, where $\iota_k\in I_\infty$ and the operators $\tilde f_\iota$ are those defined in~\ref{kashiop} together with the original ones for $\iota=i\in I^\textup{re}$.
Define the following set:\[
\mathcal B(\infty)=\{\tilde f_{\iota_1}\ldots\tilde f_{\iota_s}.1\mid\iota_k\in I_\infty\}\subset \dfrac{\mathcal L(\infty)}{v^{-1}\mathcal L(\infty)}.\]
The following theorem will be proved in section~\ref{grandloop}:
\theo\label{criinfini} The Kashiwara operators are still defined on $\mathcal B(\infty)$, which is a crystal once equipped with the following maps:\[
\textup{wt}(b)&=\sum_{i\in I}\nu_i\aaa_i\text{ if }|b|=\nu\in-\NN I\\
\eee_i(b)&=\overline{\max\{\bar{\textup{\sffamily c}}\mid \tilde e_{i,{\textup{\sffamily c}}}(b)\neq0\}}.\] 
\etheo

We have the following characterization, analogous to~\cite[3.2.3]{kashisaito}:

\prop\label{charackashisaito}
Let $\mathcal B$ be a crystal, and $b_0\in\mathcal B$ with weight $0$, such that:\benu
\item$\textup{wt}(\mathcal B)\subset-\sum_{i\in I}\NN\aaa_i$;
\item the only element of $\mathcal B$ with weight $0$ is $b_0$;
\item $\eee_i(b_0)=0$ for every $i\in I$;
\item there exists a strict embedding $\Psi_i:\mathcal B\rarrow\mathcal B_i\otimes\mathcal B$ for every $i\in I$;
\item for every $b\neq b_0$ there exists $i\in I$ such that $\Psi_i(b)=\textup{\sffamily c}\otimes b'$ for some $b'\in\mathcal B$ and $\textup{\sffamily c}\in\mathcal C_i\setminus\{(0)_i\}$.
\item for every $i$, the crystal $\mathcal B'_i= \pi_i\Psi_i(\mathcal B)$ is normal, where $\pi_i$ is the second projection $\mathcal B_i\otimes\mathcal B\rarrow\mathcal B$. \eenu
Then $\mathcal B\simeq\mathcal B(\infty)$.
\eprop

\rema The crystal structure we consider on $\mathcal B'_i$ is $(\pi_i\Psi_i)_*\mathcal B$. If $b'\in\mathcal B'_i$, we get $\tilde e_\iota(b')=0$ (respectively $\tilde f_\iota(b')=0$) if, with respect to the structure of $\mathcal B$, $\tilde e_\iota(b')\in\mathcal B\setminus\mathcal B'_i$ (respectively $\tilde f_\iota(b')\in\mathcal B\setminus\mathcal B'_i$).
\erema

\demo
First note that we necessarily have $\Psi_i(b_0)=(0)_i\otimes b_0$, thanks to (1).
Let us show that that for any $b\in \mathcal B\setminus\{b_0\}$ there exists $\iota\in I_\infty$ such that $\tilde e_\iota(b)\neq0$. Consider $i\in I$ such that $\Psi_i(b)=\textup{\sffamily c}\otimes b'$ for some $b'\in\mathcal B$ and nontrivial $\textup{\sffamily c}\in\mathcal C_i$, and assume $i$ is imaginary since the result is already known from~\cite[3.2.3]{kashisaito} when $i$ is real.
If $b'=b_0$, since $\phi_i(b_0)=0$, we have $\tilde e_{i,\textup{\sffamily c}_1}(b)=
\textup{\sffamily c}\backslash\textup{\sffamily c}_1\otimes b_0\neq0$. Otherwise  by induction on the weight, we can assume that there exists $\iota\in I_\infty$ such that $\tilde e_\iota(b')\neq0$. If $\iota=(j,1)$ for some real vertex $j$, we get $\tilde e_\iota(b)\neq0$. If $\iota=(j,l)$ for some imaginary vertex $j$, we have to prove that $\phi_j(b')=+\infty$ to get to the same result. But we have $b'=\tilde f_{j,l}\tilde e_{j,l}b'\neq0$, hence $\phi_j(\tilde e_{j,l}(b'))\neq0$ by normality of $\mathcal B'_j$. Since $(j,j)≤0$, we have:\[
\text{wt}_j(b')=\text{wt}_j(\tilde e_{j,l}(b'))-l\langle e_j,\aaa_j\rangle≥\text{wt}_j(\tilde e_{j,l}(b'))>0,\]
hence $\phi_j(b')=+\infty$, and:\[
\Psi_i(\tilde e_{j,l}(b))=\tilde e_{j,l}(\textup{\sffamily c}\otimes b')=\textup{\sffamily c}\otimes\tilde e_{j,l}(b')\neq0,\]
which proves that $\tilde e_{j,l}(b)\neq0$.

Hence any element can be written $b=\tilde f_{\iota_1}\ldots\tilde f_{\iota_1}(b_0)$ for some $\iota_k\in I_\infty$. The end of the proof is analogous to the one given in~\cite{kashisaito}, one just has to replace $I$ by $I_\infty$ (which is countably infinite).
\edemo

\subsubsection{Geometric realization}

\nota\label{notacrysgeom}
From~\ref{generA}, we have the following bijections:
\[\xymatrix{\Irr\Lambda(\nu)_{i,l}\ar[rr]_-{\sim}^-{\mathfrak k_{i,l}}&&\Irr\Lambda(\nu-le_i)_{i,0}\times\mathcal C_{i,l}}\]
where $\nu\in P^+=\NN^I$, $i\in I$, $l>0$. 
Set, for $\textup{\sffamily c}\in\mathcal C_{i,l}$: \begin{align*}
\Irr\Lambda_{i,l}&=\bigsqcup_{\nu\in P^+}\Irr\Lambda(\nu)_{i,l}\\
\Irr\Lambda(\nu)_{i,\textup{\sffamily c}}&=\mathfrak l_{i,l}^{-1}(\Irr\Lambda(\nu-le_i)_{i,0}\times\{\textup{\sffamily c}\})\\
\Irr\Lambda_{i,\textup{\sffamily c}}&=\bigsqcup_{\nu\in P^+}\Irr\Lambda(\nu)_{i,\textup{\sffamily c}}\\
\Irr\Lambda&=\bigsqcup_{\nu\in P^+}\Irr\Lambda(\nu)\end{align*}
and denote by $\tilde{ e}_{i,\textup{\sffamily c}}$ and $\tilde{ f}_{i,\textup{\sffamily c}}$ the inverse bijections:
\[\xymatrix{\tilde{ e}_{i,\textup{\sffamily c}}:\Irr\Lambda_{i,\textup{\sffamily c}}\ar@<.7ex>[rr]&&\Irr\Lambda_{i,0}:\tilde{ f}_{i,\textup{\sffamily c}}\ar@<.7ex>[ll]}\]
induced by $\mathfrak k_{i,l}$. Then, for every $l>0$, we define:\begin{align*}
\tilde{ e}_{i,l}&=\bigsqcup_{\textup{\sffamily c}\in\mathcal C_i}\tilde{ f}_{i,\textup{\sffamily c}\setminus l}\tilde{ e}_{i,\textup{\sffamily c}}:\Irr\Lambda\rarrow\Irr\Lambda\sqcup\{0\}\\
\tilde{ f}_{i,l}&=\tilde{ f}_{i,(l)}\sqcup\Big(\bigsqcup_{\textup{\sffamily c}\in\mathcal C_{i}}\tilde{ f}_{i,(l,\textup{\sffamily c})}\tilde{ e}_{i,\textup{\sffamily c}}\Big):\Irr\Lambda\rarrow \Irr\Lambda\sqcup\{0\}
\end{align*}
where $\tilde{ f}_{i,\textup{\sffamily c}\setminus l}=0$ if $\omega_i≥2$ and $l\neq \textup{\sffamily c}_1$, or if $\omega_i=1$ and $m_l( \textup{\sffamily c})=0$.
\enota

It is obvious from the definitions that we have:

\prop The set $\Irr\Lambda$ is a crystal with respect to $\textup{wt}:Z\in\Irr\Lambda(\nu)\mapsto-C\nu$, $\eee_i$ the composition of $\sqcup_{l>0}\mathfrak k_{i,l}$ and the second projection, and $\tilde{ e}_{i,l},\tilde{ f}_{i,l}$ the maps defined above.
\eprop

The duality $\Lambda\rarrow\Lambda$, $x\mapsto x^*$ induces a bijection $*:\Irr\Lambda\rarrow\Irr\Lambda$, $Z\mapsto Z^*$ preserving the grading.
Following~\cite{kashisaito}, we note:\[
\eee_i^*&=*\eee_i*\\
\tilde e_{i,l}^*&=*\tilde e_{i,l}*\\
\tilde f_{i,l}^*&=*\tilde f_{i,l}*.\]
Note that $\eee_i^*(Z)$ is the dimension of the largest subspace of $\cap_{h\in H_i}\ker x_h$ stable by $(x_h)_{h\in H_i}$, for a generic element $x\in Z$. We will note $\tilde e_{i}^{*\textup{max}}(Z)$ instead of $\tilde e_{i,\textup{\sffamily c}}(Z)$ when $\textup{\sffamily c}=\eee_i^*(Z)$. We have the following, corresponding to~\cite[5.3.1]{kashisaito} when $i$ is real:

\prop\label{geomkashisaito} Consider $Z\in \Irr\Lambda(\nu)$ such that $\eee_i^*(Z)=\textup{\sffamily c}\neq0$ for some imaginary vertex $i$, and set $\bar Z=\tilde e_{i,\textup{\sffamily c}}^*(Z)$. Assume $\textup{wt}_i(\bar Z)>0$.
We have:\benu\itemsep0.5em
\item $\eee_i(Z)=\eee_i(\bar Z)$;
\item $\left\{\begin{aligned} & \eee_i^*(\tilde e_\iota(Z))=\eee_i^*(Z) \\
&  \tilde e_{i}^{*\textup{max}}(\tilde e_\iota(Z))=\tilde e_\iota(\bar Z)\end{aligned}\right.$ for every $\iota\in I_\infty$.
\eenu
\eprop

\demo The proof is actually simpler than in the real case. Indeed, in the proof of~\ref{generA}, consider $y\in\bar Z^*$ and $z\in\textup{\sffamily c}^*$ (we abusively identify $\Irr\Lambda(le)$ with $\mathcal C_{i,l}$). We want:\[
0=\dsum_{h\in H_i}y_{\bar h}\eta_h+\dsum_{h\in\Omega(i)}\big[(y_{\bar h}\eta_h-\eta_hz_{\bar h})+(y_{h}\eta_{\bar h}-\eta_{\bar h}z_{h})\big].\]
Note that:\[
0<\textup{wt}_i(\bar Z)=\sum_{h\in H_i}\nu_{t(h)}-(i,i)(\nu_i-|\textup{\sffamily c}|)\Lrarrow0<\sum_{h\in H_i}\nu_{t(h)}\]
since $(i,i)≤0$, and since it is impossible to have $\sum_{h\in H_i}\nu_{t(h)}=0$ and $\nu_i-|\textup{\sffamily c}|>0$. Hence there exists $h_0\in H_i$ such that $\nu_{h_0}>0$.
We have $\Spec(z_{\bar h_1})\cap\Spec(y_{\bar h_1})=\varnothing$ for a generic choice of $y,z$, where $h_1\in\Omega(i)$. But then the map:\[
\eta_{h_1}\mapsto y_{\bar h_1}\eta_{h_1}-\eta_{h_1}z_{\bar h_1}\]
is invertible, and we can generically choose $\eta_{h_0}$ so that:\[
\dim\CC\langle z^*_h\mid  h\in H_i\rangle.\ima \eta^*_{h_0}=|\textup{\sffamily c}|.\]
 This proves $\eee_i^*(Z^*)=\eee_i^*(\bar Z^*)$, hence (1). The second statement directly follows from the proof of (1).
\edemo

\theo\label{isocrys} We have $\Irr\Lambda\simeq\mathcal B(\infty)$.\etheo

\demo Set $\Psi_i(Z)=\eee_i^*(Z)\otimes\tilde e_i^{*\textup{max}}(Z)$, which is clearly injective. By~\ref{geomkashisaito} and the definition of our generalized crystals, $\Psi_i$ is a morphism. Note that we have $\Psi(\tilde e_{i,l}(Z))=\tilde e_{i,l}(\eee^*_i(Z))\otimes\bar Z$ if $\textup{wt}_i(\bar Z)=0$. By~\ref{remlu} (or its dual analog), the condition (5) of~\ref{charackashisaito} is satisfied. The condition (6) is satisfied because it is clear that $\tilde f_{i,l}(\bar Z)\notin\mathcal B'_i$ if $\phi_i(\bar Z)=0$. Hence we get the result.
\edemo

\subsubsection{Semicanonical basis}

The following proposition is proved in~\cite{article2}:

\prop\label{surjsemican} There exists a surjective morphism $\Phi:U^+_{v=1}\rarrow\mathcal M_\circ$ defined by:\[
\left\{\begin{aligned}& E_{i,a}\mapsto 1_{i,l}  &&   \text{if }i\in I^\textup{im}     \\  
&E_{i}\mapsto 1_i &&   \text{if }i\in I^\textup{re}. \end{aligned}\right.\]
\eprop

 Thanks to~\ref{isocrys}, we now have:
 \theo The morphism $\Phi$ is an isomorphism $U^+_{v=1}\isom\mathcal M_\circ$.\etheo
 
 \demo The family $(f_Z)_{Z\in\Irr\Lambda}$ is clearly free, so we have:\[
 |\Irr\Lambda(\nu)|≤\dim\mathcal M_\circ(\nu)≤\dim U^+_{v=1}[\nu], \]
 the latter inequality being true thanks to~\ref{surjsemican}. From~\ref{isocrys}, we have $|\Irr\Lambda(\nu)|=\dim U^+_{v=1}[\nu]$, hence $(f_Z)_{Z\in\Irr\Lambda}$ is a basis of $\mathcal M_\circ$, and $\Phi$ is an isomorphism.
 \edemo

\defi The semicanonical basis of $U^+_{v=1}$ is the pullback of $(f_Z)_{Z\in\Irr\Lambda}$.\edefi

\subsection{The crystals $\mathcal B(\lambda)$}

\subsubsection{Algebraic definition}

We will use the \emph{fundamental weights} $(\Lambda_i)_{i\in I}$ defined by $(i,\Lambda_j)=\delta_{i,j}$ for every $i,j\in I$. Note that the isomorphism $P\isom\sum\ZZ\Lambda_i$, $e_i\mapsto\Lambda_i$ maps $\aaa_i$ to $i$. We use this isomorphism to identify $\sum\ZZ\Lambda_i$ with $P$ and $\sum\NN\Lambda_i$ with $P^+$. We call \emph{dominant} the elements $\lambda\in P^+$, which are the ones satisfying $(i,\lambda)≥0$ for every $i\in I$.

\defi 
We denote by $\mathcal O$ the category of $U$-modules satisfying:\benu\itemsep0.5em
\item $M=\oplus_{\mu\in P}M_\mu$ where $M_\mu=\big\{m\in M\mid \forall i,~K_im=v^{(\mu,i)}m\big\}$;
\item For any $m\in M$, there exists $p≥0$ such that $xm=0$ as soon as $x\in U^+[\nu]$ and $\text{ht}(\nu)≥p$.\eenu
 For any $\lambda\in P$, we define a \emph{Verma module}:\[
M(\lambda)=\dfrac{U}{\dsum_{\iota\in I_\infty}UE_\iota+\dsum_{i\in I}U\big(K_i-v^{(i,\lambda)}\big)}\in\mathcal O\]
and  the following simple quotient:\[
\pi_\lambda:U\twoheadrightarrow V(\lambda)=\dfrac{M(\lambda)}{M(\lambda)^-}\in\mathcal O\]
where $M(\lambda)^-$ is the sum of all strict submodules of $M(\lambda)$. We will denote by $v_\lambda\in V(\lambda)_\lambda$ the image of $1\in U$.
\edefi
\rema\label{remtriang} Note that thanks to~\ref{kashiprop} (2), we have a triangular decomposition, and thus $M(\lambda)=U^-v_\lambda$.\erema

We have the following proposition, generalizing the case $i\in I^\text{re}$:

\prop\label{kashioplambda} Assume $(i,\lambda)≥0$ for some imaginary vertex $i$. Then we have the following decomposition:\[
V(\lambda)=\bigoplus_{\textup{\sffamily c}\in\mathcal C_i}b_{i,\textup{\sffamily c}}\mathcal K_i\]
where $\mathcal K_i=\dcap_{l>0}\ker E_{i,l}$.
\eprop

\demo
Let's first prove the existence. Consider $v\in V(\lambda)$, and assume first that $v$ is of the following form: $v=ub_{i,\textup{\sffamily c}}z$ for some $\textup{\sffamily c}\in\mathcal C_i$, $u\in U^-$ satisfying $[a_{i,l},u]=0$ for every $l$, and $z\in \mathcal K_i$. We proceed by induction on $|\textup{\sffamily c}|$. 
First note that if $(i,|u|)=0$, since $i$ is imaginary one necessarily has:\[
\text{supp}|u|\subseteq\{j\in I\mid (i,j)=0\}.\]
Hence $[u,b_{i,l}]=0$ for any $l$ (whether $i$ is isotropic or not) and we get the result by induction.
Otherwise, $(i,|u|)>0$, and we set:\[
l&=\textup{\sffamily c}_1\\
[u,b_{i,l}]^\circ&=ub_{i,l}-R(v)b_{i,l}u\text{ for some }R\in\QQ(v)\\
z'&=b_{i,\textup{\sffamily c}\backslash\textup{\sffamily c}_1}z\in V(\lambda)_\mu.\] 
For any $k>0$, we have:\[
[a_{i,k},[u,b_{i,l}]^\circ]z'&=\delta_{l,k}\tau_{i,l}\big\{u(K_{-li}-K_{li})-R(v)(K_{-li}-K_{li})u\big\}z'\\
&\qquad[~\textit{c.f.}~\ref{kashiprop}~(2)~]\\
&=\delta_{l,k}\tau_{i,l}u\big\{(v^{-l(i,\mu)}-v^{l(i,\mu)})\\
&\qquad\qquad\qquad-R(v)(v^{-l(i,|u|+\mu)}-v^{l(i,|u|+\mu)})\big\}z'\\
&=0\]
if:\[
R(v)=\dfrac{v^{-l(i,\mu)}-v^{l(i,\mu)}}{v^{-l(i,|u|+\mu)}-v^{l(i,|u|+\mu)}},\]
which is possible since:\[
(i,|u|+\mu)=(i,\lambda)+(i,|u|)+(i,\mu-\lambda)>(i,\lambda)+(i,\mu-\lambda)≥0.\]
We have used that since $i$ is imaginary, we have:\[
\mu-\lambda\in-\NN I\Rarrow(i,\mu-\lambda)≥0.\]
Hence, the following equality:\[
v=[u,b_{i,l}]^\circ b_{i,\textup{\sffamily c}\backslash\textup{\sffamily c}_1}z+R(v)b_{i,l}ub_{i,\textup{\sffamily c}\backslash\textup{\sffamily c}_1}z
\]
along with the induction hypothesis allow us to conclude since $|\textup{\sffamily c}\backslash\textup{\sffamily c}_1|<|\textup{\sffamily c}|$, and since $\oplus_{\textup{\sffamily c}\in\mathcal C_i}b_{i,\textup{\sffamily c}}\mathcal K_i$ is stable by left-multiplication by $b_{i,l}$.

Then, we prove the existence of the decompostion for a general $v\in V(\lambda)_\mu$, using induction on $\sum(\lambda_i-\mu_i)$. If $v\neq v_\lambda$, thanks to~\ref{remtriang}, we can write:\[
v=\dsum_{\iota\in I_\infty}b_{\iota}v_\iota\]
for some finitely many nonzero $v_\iota\in V(\lambda)$. Thanks to our induction hypothesis, we have:\[
v=\dsum_{\iota\in I_\infty,\textup{\sffamily c}\in\mathcal C_i}b_\iota b_{i,\textup{\sffamily c}} z_{\iota,\textup{\sffamily c}}\]
for some finitely many nonzero $z_{\iota,\textup{\sffamily c}}\in\mathcal K_i$. Then:\[
v=\dsum_{l>0,\textup{\sffamily c}\in\mathcal C_i}b_{i,(l,\textup{\sffamily c})}z_{(i,l),\textup{\sffamily c}}+\dsum_{\substack{\iota\in I_\infty\backslash(\{i\}\times\NN_{>0})\\\textup{\sffamily c}\in\mathcal C_i}}b_\iota b_{i,\textup{\sffamily c}} z_{\iota,\textup{\sffamily c}}\]
and we have the result since $b_\iota b_{i,\textup{\sffamily c}} z_{\iota,\textup{\sffamily c}}$ is of the form $ub_{i,\textup{\sffamily c}}z$ already treated. Indeed, thanks to~\ref{kashiprop} (2), $[a_{i,l},b_{j,k}]=0$ for any $l,k>0$ if $j\neq i$.

To prove the unicity of the decomposition, consider a minimal nontrivial relation of dependance:\[
0=\dsum_{\textup{\sffamily c}\in\mathcal C_i}b_{i,\textup{\sffamily c}}z_\textup{\sffamily c},\]
where $z_\textup{\sffamily c}\in V(\lambda)_{\mu+|\textup{\sffamily c}|i}\cap\mathcal K_i$. 
We have to considerate separately the cases $i\in I^\textup{iso}$ and $i\notin I^\textup{iso}$.
First, consider $i\notin I^\textup{iso}$.
Consider $r$ maximal such that there exists $\textup{\sffamily c}=(\textup{\sffamily c}_1,\ldots,\textup{\sffamily c}_r)$ such that $z_\textup{\sffamily c}\neq0$.
Set for any $k\in[\![1,r]\!]$:\[
\textup{\sffamily c}_{<k}&=(\textup{\sffamily c}_1,\ldots,\textup{\sffamily c}_{k-1})\\
\textup{\sffamily c}_{>k}&=(\textup{\sffamily c}_{k+1},\ldots,\textup{\sffamily c}_{r})\]
with the convention $\textup{\sffamily c}_{<1}=\varnothing=\textup{\sffamily c}_{>r}$. 
Then, if $l>0$, we get the following from~\ref{kashiprop} (2), where by convention $b_{i,\varnothing}=1$:\[
[a_{i,l},b_{i,\textup{\sffamily c}}]&=\tau_{i,l}\dsum_{k:\textup{\sffamily c}_k=l}b_{i,\textup{\sffamily c}_{<k}}(K_{-li}-K_{li})b_{i,\textup{\sffamily c}_{>k}}\\
&=\tau_{i,l}\dsum_{k:\textup{\sffamily c}_k=l}b_{i,\textup{\sffamily c}\backslash\textup{\sffamily c}_{k}}(v_i^{2l|\textup{\sffamily c}_{>k}|}K_{-li}-v_i^{-2l|\textup{\sffamily c}_{>k}|}K_{li}).
\]
Then, since $z_\textup{\sffamily c}\in\mathcal K_i$:\[
a_{i,l}b_{i,\textup{\sffamily c}}z_\textup{\sffamily c}&=\tau_{i,l}\dsum_{k:\textup{\sffamily c}_k=l}b_{i,\textup{\sffamily c}\backslash\textup{\sffamily c}_{k}}(v_i^{2l|\textup{\sffamily c}_{>k}|}K_{-li}-v_i^{-2l|\textup{\sffamily c}_{>k}|}K_{li})z_\textup{\sffamily c}\\
&=\tau_{i,l}\dsum_{k:\textup{\sffamily c}_k=l}b_{i,\textup{\sffamily c}\backslash\textup{\sffamily c}_{k}}(v^{-l(i,\mu+|\textup{\sffamily c}_{≤k}|i)}-v^{l(i,\mu+|\textup{\sffamily c}_{≤k}|i)})z_\textup{\sffamily c}
\]
where $\textup{\sffamily c}_{≤k}=(\textup{\sffamily c}_{<k},\textup{\sffamily c}_{k})$.
We see that for any $\textup{\sffamily c}'\in\mathfrak S_r\textup{\sffamily c}$ (with the convention $(\sigma\textup{\sffamily c})_k=\textup{\sffamily c}_{\sigma(k)}$), since $r$ is maximal, we have:\[
0=a_{i,\textup{\sffamily c}'}\dsum_{\textup{\sffamily c}''\in\mathcal C_i}b_{i,\textup{\sffamily c}''}z_{\textup{\sffamily c}''}&=a_{i,\textup{\sffamily c}'}\dsum_{\textup{\sffamily c}''\in\mathfrak S_r\textup{\sffamily c}}b_{i,\textup{\sffamily c}''}z_{\textup{\sffamily c}''}\\
&=\tau_{i,\textup{\sffamily c}}\dsum_{\textup{\sffamily c}''\in\mathfrak S_r\textup{\sffamily c}}P_{\textup{\sffamily c}',\textup{\sffamily c}''}(v)z_{\textup{\sffamily c}''}\]
where $P_{\textup{\sffamily c}',\textup{\sffamily c}''}(v)\in\ZZ[v,v^{-1}]$. Since $(z_{\textup{\sffamily c}''})_{\textup{\sffamily c}''\in\mathfrak S_r\textup{\sffamily c}}\neq0$, we have to prove that:\[
\Delta(v)=\det(P_{\textup{\sffamily c}',\textup{\sffamily c}''}(v))_{\textup{\sffamily c}',\textup{\sffamily c}''\in\mathfrak S_r\textup{\sffamily c}}\neq 0\in\ZZ[v,v^{-1}]\]
to end our proof in the case $(i,i)<0$.
Note that $\lambda-(\mu+|\textup{\sffamily c}|i)\in\NN I$, hence, since $i$ is imaginary:\[
(i,\mu+|\textup{\sffamily c}|i)=(i,\lambda)+(i,\mu+|\textup{\sffamily c}|i-\lambda)≥0.\]
Then, for any $\textup{\sffamily c}'\in\mathfrak S_r\textup{\sffamily c}$, one has:\[
\underset{\textup{\sffamily c}''\in\mathfrak S_r\textup{\sffamily c}}{\textup{max}}\{\textup{deg}(P_{\textup{\sffamily c}',\textup{\sffamily c}''})\}=\textup{deg}(P_{\textup{\sffamily c}',\bar{\textup{\sffamily c}}'})=\dsum_{1≤k≤r}\textup{\sffamily c}_k(i,\mu+\textup{\sffamily c}_ki)= m\]
which is only reached for $\textup{\sffamily c}''=\bar{\textup{\sffamily c}}'$. However, this is not true if $m=0$, which can only happen if our initial relation of dependance is of the form $b_{i,l}z_{l}=0$, with $(i,\mu+li)=0$. But if $(i,\mu+li)=0$, the module generated by $b_{i,l}z_l$ is a nontrivial strict submodule of $V(\lambda)$ since for every $k>0$ and $j\neq i$:\[
a_{i,k} b_{i,l}z_l&=0\\
 a_{j,k} b_{i,l}z_l&=b_{i,l}a_{j,k}z_l.\] 
 Hence the relation of dependance $b_{i,l}z_{l}=0$ is actually trivial by definition of $V(\lambda)$.

Otherwise, the application $\mathfrak S_r\textup{\sffamily c}\rarrow\mathfrak S_r\textup{\sffamily c}$, ${\textup{\sffamily c}}'\mapsto\bar{\textup{\sffamily c}}'$ being a permutation, the degree of $\Delta$ is $|\mathfrak S_r\textup{\sffamily c}|m>0$, and in particular $\Delta\neq0$.

We finally have to prove the uniqueness in the case $(i,i)=0$. Write a relation of dependance of minimal degree:\[
0=\dsum_{\nu\in\mathcal C_i}b_{i,\nu}z_\nu,\]
where $z_\nu\in V(\lambda)_{\mu+|\nu|i}\cap\mathcal K_i$. For any $l>0$, we have, thanks to~\ref{kashiprop} (2):\[
a_{i,l}\dsum_{\nu\in\mathcal C_i}b_{i,\nu}z_\nu&=\dsum_{\nu\in\mathcal C_i}m_l(\nu)\tau_{i,l}\big(v^{-l(i,\mu)}-v^{l(i,\mu)}\big)b_{i,\nu\backslash l}z_\nu\\
&=\tau_{i,l}\big(v^{-l(i,\mu)}-v^{l(i,\mu)}\big)\dsum_{\nu'\in\mathcal C_i}b_{i,\nu'}\big\{(m_l(\nu')+1)z_{\nu'\cup l}\big\}\]
which contradicts the minimality of the first relation. Note that we can assume that $l(i,\mu)\neq0$: otherwise we would again have an initial trivial relation of dependance (more precisely for every $\nu$ we would have $b_{i,\nu}z_\nu=0\in V(\lambda)$).
\edemo

This proposition allows us to define Kashiwara operators $\tilde e_\iota$, $\tilde f_\iota$ on each $V(\lambda)$, exactly as in~\ref{kashiop}:

\defi\label{kashiopbis}
If $i$ is imaginary and $v=\sum_{\textup{\sffamily c}\in\mathcal C_i}b_{i,\textup{\sffamily c}}z_\textup{\sffamily c}\in V(\lambda)$, set:\[
\tilde e_{i,l}(v)&=\left\{\begin{aligned}&\dsum_{\textup{\sffamily c}:\textup{\sffamily c}_1=l}b_{i,\textup{\sffamily c}\backslash\textup{\sffamily c}_1}z_\textup{\sffamily c}&&\text{if }i\notin I^\text{iso}\\
&\dsum_{\nu\in\mathcal C_i}\sqrt{\dfrac{m_l(\nu)}{l}}b_{i,\nu\backslash l}z_\nu&&\text{if }i\in I^\text{iso}\end{aligned}\right.\\
\tilde f_{i,l}(v)&=\left\{\begin{aligned}&\dsum_{\textup{\sffamily c}\in\mathcal C_i}b_{i,(l,\textup{\sffamily c})}z_\textup{\sffamily c}&&\text{if }i\notin I^\text{iso}\\
&\dsum_{\nu\in\mathcal C_i}\sqrt{\frac{l}{m_l(\nu)+1}}b_{i,\nu\cup l}z_\nu&&\text{if }i\in I^\text{iso}.\end{aligned}\right.\]
\edefi

The following will be proved in section~\ref{grandloop}:

\theo\label{crilambda} Assume $\lambda$ is dominant. The Kashiwara operators, along with the maps:\[
\textup{wt}(m)&=\mu\text{ if }m\in V(\lambda)_\mu\\
\eee_i(m)&=\overline{\max\{\bar{\textup{\sffamily c}}\mid \tilde e_{i,{\textup{\sffamily c}}}(m)\neq0\}},\]
induce a structure of crystal on:\[
\mathcal B(\lambda)=\{\tilde f_{\iota_1}\ldots\tilde f_{\iota_s}v_\lambda\mid\iota_k\in I_\infty\}\subset\dfrac{\mathcal L(\lambda)}{v^{-1}\mathcal L(\lambda)}\]
where:\[
\mathcal L(\lambda)=\dsum_{\iota_1,\ldots,\iota_s}\mathcal A\tilde f_{\iota_1}\ldots\tilde f_{\iota_s}v_\lambda.\]
\etheo

\rema These crystals are normal: consider $m\in V(\lambda)_\mu$ and $i$ imaginary (again, the case of real vertices is already known). We have already seen that we necessarily have $(i,\mu)≥0$ since $\lambda$ is dominant. If $(i,\mu)=0$ for some imaginary vertex $i$, then $a_\iota b_{i,l}m=b_{i,l}a_\iota m$ for any $\iota\in I_\infty$ (use~\ref{kashiprop} (2) if $\iota=(i,l)$). Hence, for any $l>0$, the submodule of $V(\lambda)$ spanned by $b_{i,l}m$ is a strict submodule, and we get $\tilde f_{i,l}m=0$.

Otherwise, $(i,\mu)>0$, and for every $\mu'\in-\NN I$, since $(i,i)≤0$, we get:\[
(i,\mu+\mu')=(i,\mu)+(i,\mu')≥(i,\mu)>0,\]
hence $\text{max}\{|\text{\sffamily c}|\mid b_{i,\text{\sffamily c}}m\neq0\}=+\infty$.
\erema

\subsubsection{Geometric realization}

\nota
Consider $\lambda$ dominant.
We have the following bijections:
\[\xymatrix{\Irr\mathfrak L(\nu,\lambda)_{i,l}\ar[rr]_-{\sim}^-{\mathfrak l_{i,l}}&&\Irr\mathfrak L(\nu-le_i,\lambda)_{i,0}\times\mathcal C_{i,l}}\]
each time the left hand side is non-empty (\textit{c.f.}~\ref{subcryst}). 
Set, for $\textup{\sffamily c}\in\mathcal C_{i,l}$: \begin{align*}
\Irr\mathfrak L(\lambda)_{i,l}&=\bigsqcup_{\nu\in P^+}\Irr\mathfrak L(\nu,\lambda)_{i,l}\\
\Irr\mathfrak L(\nu,\lambda)_{i,\textup{\sffamily c}}&=\mathfrak l_{i,l}^{-1}(\Irr\mathfrak L(\nu-le_i,\lambda)_{i,0}\times\{\textup{\sffamily c}\})\\
\Irr\mathfrak L(\lambda)_{i,\textup{\sffamily c}}&=\bigsqcup_{\nu\in P^+}\Irr\mathfrak L(\nu,\lambda)_{i,\textup{\sffamily c}}\\
\Irr\mathfrak L(\lambda)&=\bigsqcup_{\nu\in P^+}\Irr\mathfrak L(\nu,\lambda)\end{align*}
and denote by $\tilde{ e}_{i,\textup{\sffamily c}}$ and $\tilde{ f}_{i,\textup{\sffamily c}}$ the inverse bijections:
\[\xymatrix{\tilde{ e}_{i,\textup{\sffamily c}}:\Irr\mathfrak L(\lambda)_{i,\textup{\sffamily c}}\ar@<.7ex>[rr]&&\Irr\mathfrak L(\lambda)_{i,0}:\tilde{ f}_{i,\textup{\sffamily c}}\ar@<.7ex>[ll]}\]
induced by $\mathfrak l_{i,l}$. Then, for every $l>0$, we define:\begin{align*}
\tilde{ e}_{i,l}&=\bigsqcup_{\textup{\sffamily c}\in\mathcal C_i}\tilde{ f}_{i,\textup{\sffamily c}\setminus \textup{\sffamily c}_1}\tilde{ e}_{i,\textup{\sffamily c}}:\Irr\mathfrak L(\lambda)\rarrow\Irr\mathfrak L(\lambda)\sqcup\{0\}\\
\tilde{ f}_{i,l}&=\tilde{ f}_{i,(l)}\sqcup\Big(\bigsqcup_{\textup{\sffamily c}\in\mathcal C_{i}}\tilde{ f}_{i,(l,\textup{\sffamily c})}\tilde{ e}_{i,\textup{\sffamily c}}\Big):\Irr\mathfrak L(\lambda)\rarrow \Irr\mathfrak L(\lambda)\sqcup\{0\}
\end{align*}
with the same conventions than in~\ref{notacrysgeom}.
\enota

The following is a direct consequence of~\ref{subcryst}:
\prop The set $\Irr\mathfrak L(\lambda)$ is a crystal with respect to:\[
\textup{wt}:b\in\Irr\mathfrak L(\nu,\lambda)\mapsto\lambda-C\nu,\]
 $\eee_i$ the composition of $\sqcup_{l>0}\mathfrak l_{i,l}$ and the second projection, and $\tilde{ e}_{i,l},\tilde{ f}_{i,l}$ the maps defined above. \eprop

\rema Thanks to~\ref{subcryst} and the classical case, we have, for every $i\in I$:\[
\phi_i(b)=\max\{|\textup{\sffamily c}|\in\NN\mid \tilde{ f}_{i,\textup{\sffamily c}}(b)\neq0\}.\]
Indeed, for $b\in \Irr\mathfrak L(\nu,\lambda)$, it is impossible to have $\nu_i>0$ and $\lambda_i+\sum_{h\in H_i}\nu_{t(h)}=0$, hence:\[
\lambda_i+\sum_{h\in H_i}\nu_{t(h)}>0\Lrarrow \langle e_i,\lambda-C\nu\rangle >0,\]
and $\Irr\mathfrak L(\lambda)$ is normal.
\erema

In an analogous way, one can equip $\Irr\widetilde{\mathfrak Z}$ with a structure of crystal, thanks to~\ref{bijtens}, and get:

\theo The crystal structure on $\Irr\widetilde{\mathfrak Z}$ coincides with that of the tensor product $\Irr\mathfrak L(\lambda)\otimes\Irr\mathfrak L(\lambda')$.\etheo
\demo This is essentially~\ref{tenscrysgeom}.\edemo

We will see in section~\ref{grandloop} how the previous theorem leads to the following:
\theo\label{jotens} If $\lambda$ is dominant, we have the following isomorphism of crystals $\mathcal B(\lambda)\simeq\Irr\mathfrak L(\lambda)$.\etheo

\subsection{Grand loop argument}\label{grandloop}

To prove theorems~\ref{criinfini},~\ref{crilambda} and~\ref{jotens}, one has to generalize Kashiwara's grand-loop argument to our framework (see~\cite{kashiloop}).

Instead of giving the whole grand-loop argument, we give a few lemmas that yield its generalization.

\nota When working with a $\mathcal A$-lattice $\mathcal L$, we will write $m\equiv m'$ instead of $m=m'+v^{-1}\mathcal L$ for any $m,m'\in \mathcal L$.\enota

The following result is about the tensor product:
\lemm Consider two dominant weights $\lambda$ and $\lambda'$, and $(m,m')\in \mathcal L(\lambda)\times \mathcal L(\lambda')_{\mu'}$. Then, for every imaginary vertex $i$ and $l>0$, we have:\[
b_{i,l}(m\otimes m')\equiv\left\{\begin{aligned}&m\otimes b_{i,l}m'&&\text{if }\textup{wt}_i(m')>0\\
&b_{i,l}m\otimes m'&&\text{otherwise.}\end{aligned}\right.\]
\elemm

\rema Note that since $\eee_i(m)\neq-\infty$ in this situation, this is exactly~\ref{defitensorprod} (7).\erema

\demo We have already seen that when $i\in I^\text{im}$, since $\mu'-\lambda'\in-\NN I$, we have:\[
(i,\mu')=(i,\lambda')+(i,\mu'-\lambda')≥0.\]
We have also already seen that thanks to~\ref{kashiprop} (2), if $(i,\mu')=0$, then $a_\iota b_{i,l}m'=b_{i,l}a_\iota m'$ for every $\iota\in I_\infty$. Hence $b_{i,l}m'=0$ since the module spanned by $b_{i,l}m'$ is a strict submodule of $V(\lambda')$.
Hence:\[
b_{i,l}(m\otimes m')&=b_{i,l} m\otimes K_{-li}m'+m\otimes b_{i,l}m'\\
&=v^{-l(i,\mu')}b_{i,l} m\otimes m'+m\otimes b_{i,l}m'\\
&\equiv\left\{\begin{aligned}&m\otimes b_{i,l}m'&&\text{if }\textup{wt}_i(m')>0\\
&b_{i,l}m\otimes m'&&\text{otherwise.}\end{aligned}\right.\]
\edemo

\lemm\label{congrutau} Consider $i\in I^\textup{im}$ and $l>0$. We have $\tau_{i,l}\equiv1/l$ if $i\in I^\textup{iso}$, $\tau_{i,l}\equiv1$ otherwise.\elemm

\demo First note that for any $i\in I^\text{im}$ and $l>0$, $\{E_{i,l},E_{i,l}\}\equiv1$ is required by~\ref{hypo}. 
Assume moreover that:\[
\{E_{i,l},E_{i,l}\}=\prod_{1≤k≤l}\dfrac{1}{1-v^{-k}},
\]
which is consistent with~\cite[2.29]{article2}. Then, when $i\in I^\text{iso}$, we have an isomorphism from the ring of symmetric functions $\Lambda=\ZZ[x_k,k≥1]$ to $\ZZ[E_{i,l},l≥1]$ mapping the elementary symmetric functions $e_l$ to $v^{-l/2}E_{i,l}$ and such that the pushforward of $\{-,-\}$ is the Hall-Littlewood scalar product (still denoted by $\{-,-\}$). Asking for $a_{i,l}$ to be primitive and to satisfy $E_{i,l}-a_{i,l}\in\QQ(v)[E_{i,k},k<l]$ means that the power sum symmetric functions $p_l$ are mapped to $v^{-l/2}(-1)^{l-1}la_{i,l}$. Since the Hall-Littlewood scalar product satisfy $\{p_l,p_l\}=\frac{lv^{-l}}{1-v^{-l}}$, we get:\[
\tau_{i,l}=\{a_{i,l},a_{i,l}\}=\frac{v^l}{l^2}\frac{lv^{-l}}{1-v^{-l}}=\frac{1/l}{1-v^{-l}}\equiv1/l\]
as expected.

If $i\notin I^\text{iso}$, let us prove that $a_{i,l}\equiv E_{i,l}$ by induction on $l$. Write:\[
a_{i,l}-E_{i,l}=\dsum_{\text{\sffamily c}\in\mathcal C_{i,l}\backslash\{(l)\}}\aaa_\text{\sffamily c}a_{i,\text{\sffamily c}}\]
for some $\aaa_\text{\sffamily c}\in\QQ(v)$. By~\ref{prim}, for every $\text{\sffamily c}'\in\mathcal C_{i,l}\backslash\{(l)\}$, we have:\[
\dsum_{\text{\sffamily c}\in\mathcal C_{i,l}\backslash\{(l)\}}\aaa_\text{\sffamily c}\{a_{i,\text{\sffamily c}},a_{i,\text{\sffamily c}'}\}&=-\{E_{i,l},a_{i,\text{\sffamily c}'}\}\\
&=-\{\delta(E_{i,l}),a_{i,\text{\sffamily c}_1'}\otimes a_{i,\text{\sffamily c}'\backslash\text{\sffamily c}'_1}\}\\
&=-v_i^{\text{\sffamily c}_1'|\text{\sffamily c}'\backslash\text{\sffamily c}'_1|}\{E_{i,\text{\sffamily c}'_1},a_{i,\text{\sffamily c}'_1}\}\{E_{i,|\text{\sffamily c}'\backslash\text{\sffamily c}'_1|},a_{i,\text{\sffamily c}'\backslash\text{\sffamily c}'_1}\}\\
&=-v_i^{\sum\text{\sffamily c}'_k\text{\sffamily c}'_{k+1}}\textstyle\prod\{E_{i,\text{\sffamily c}'_k},a_{i,\text{\sffamily c}'_k}\}\\
&=-v_i^{\sum\text{\sffamily c}'_k\text{\sffamily c}'_{k+1}}\textstyle\prod\tau_{i,\text{\sffamily c}'_k}\\
&\equiv0\]
by the induction hypothesis, and since $(i,i)<0$. We have also used that $\tau_{i,k}=\{E_{i,k},a_{i,k}\}$ since $\{a_{i,k},a_{i,k}-E_{i,k}\}=0$. Also, note that:\[
\det(\{a_{i,\text{\sffamily c}},a_{i,\text{\sffamily c}'}\})_{\text{\sffamily c},\text{\sffamily c}'\in\mathcal C_{i,l}\backslash\{(l)\}}\equiv1\]
since:\[
\{a_{i,\text{\sffamily c}},a_{i,\text{\sffamily c}'}\}&=\{a_{i,\text{\sffamily c}_1}\otimes a_{i,\text{\sffamily c}\backslash\text{\sffamily c}_1},\textstyle\prod(a_{i,\text{\sffamily c}_k'}\otimes1+1\otimes a_{i,\text{\sffamily c}'_k})\}\\
&=\sum_{k:\text{\sffamily c}'_k=\text{\sffamily c}_1}v_i^{\text{\sffamily c}'_{k-1}\text{\sffamily c}'_k}\tau_{i,\text{\sffamily c}_1}\{a_{i,\text{\sffamily c}\backslash \text{\sffamily c}_1},a_{i,\text{\sffamily c}'\backslash\text{\sffamily c}'_k}\}\\
&\equiv\delta_{\text{\sffamily c}_1,\text{\sffamily c}'_1}\{a_{i,\text{\sffamily c}\backslash \text{\sffamily c}_1},a_{i,\text{\sffamily c}'\backslash\text{\sffamily c}'_k}\}\\
&\equiv\dots\equiv\delta_{\text{\sffamily c},\text{\sffamily c}'}. 
\]
Hence we get $\aaa_\text{\sffamily c}\equiv0$, which implies $\tau_{i,l}=\{a_{i,l},E_{i,l}\}\equiv\{E_{i,l},E_{i,l}\}\equiv1$.
\edemo

The following lemma deals with the behaviour of the generalized Kashiwara operators regarding our Hopf bilinear form $\{-,-\}$:

\lemm For any $u,v\in U^-$ and $(i,l)\in I_\infty$, $\{\tilde f_{i,l}u,v\}\equiv\{u,\tilde e_{i,l}v\}$.\elemm 

\demo
We can assume that $v=b_{i,\textup{\sffamily c}}z$ for some $z\in\mathcal K_i$. If $i\in I^\text{im}$ is not isotropic:\[
\{\tilde f_{i,l}u,v\}&=\{ b_{i,l}u,v\}\\
&=\{b_{i,l}\otimes u,\delta(b_{i,\textup{\sffamily c}}z)\}\\
&=\tau_{i,l}\{u,\delta^{i,l}(b_{i,\textup{\sffamily c}}z)\}\\
&=\tau_{i,l}\{u,\delta^{i,l}(b_{i,\textup{\sffamily c}})z\}\qquad[\textit{c.f.~}\text{\ref{kashiprop} (1)}]\\
&\equiv\{u,\delta^{i,l}(b_{i,\textup{\sffamily c}})z\}\qquad[\textit{c.f.~}\ref{congrutau}]\\
&=\sum_{k:\textup{\sffamily c}_k=l}v_i^{l\textup{\sffamily c}_{k-1}}\{u,b_{i,\textup{\sffamily c}\backslash\textup{\sffamily c}_k}z\}\qquad[\textit{c.f.~}\text{\ref{kashiprop} (3)}]\\
&\equiv\{u,b_{i,\textup{\sffamily c}\backslash\textup{\sffamily c}_1}z\}\\
&=\{u,\tilde e_{i,l}v\}.\]
This computation also proves the case $i\in I^\text{re}$, $l=1$, which is already known.
If $i$ is isotropic, $v=b_{i,\nu}z$ and $u=b_{i,\nu'}z'$:\[
\{\tilde f_{i,l}u,v\}&=\sqrt\frac{l}{m_l(\nu')+1}\{ b_{i,l}u,v\}\\
&=\sqrt\frac{l}{m_l(\nu')+1}\{b_{i,l}\otimes u,\delta(b_{i,\nu}z)\}\\
&=\sqrt\frac{l}{m_l(\nu')+1}\tau_{i,l}\{u,\delta^{i,l}(b_{i,\nu}z)\}\\
&=\sqrt\frac{l}{m_l(\nu')+1}\tau_{i,l}\{u,\delta^{i,l}(b_{i,\nu})z\}\\
&=\sqrt\frac{l}{m_l(\nu')+1}\tau_{i,l}m_l(\nu)\{u,b_{i,\nu\backslash l}z\}.\]
We see by induction that $\{\tilde f_{i,l}u,v\}=\{u,\tilde e_{i,l}v\}=0$ if $\nu\neq\nu'\cup l$. Otherwise, we get, thanks to~\ref{congrutau}:\[
\{\tilde f_{i,l}u,v\}\equiv\sqrt\frac{m_l(\nu)}{l}\{u,b_{i,\nu\backslash l}z\}=\{u,\tilde e_{i,l}v\}.\]
\edemo

In order to get an analogous result regarding the lattices $\mathcal L(\lambda)$, first note that there exists for each $\lambda\in P^+$ a unique symmetric bilinear form $(-,-)$ on $V(\lambda)$ satisfying:\[
(K_iu,u')&=(u,K_iu')\\
(b_{i,l}u,u')&=-(u,K_{-li}a_{i,l}u')\text{ if }i\in I^\text{im}\\
(b_iu,u')&=\dfrac{v}{v^{-1}-v}(u,K_{-i}a_{i}u')\text{ if }i\in I^\text{re}\\
(v_\lambda,v_\lambda)&=1\]
for every $u,u'\in V(\lambda)$ and $(i,l)\in I_\infty$. Then:

\lemm For every $u,v\in\mathcal L(\lambda)$ and $(i,l)\in I_\infty$, $(\tilde f_{i,l}u,v)\equiv(u,\tilde e_{i,l}v)$.\elemm

\demo Assume that $v=b_{i,\text{\sffamily c}}z$ for some nontrivial $\text{\sffamily c}$ and some $z\in\mathcal K_i\cap V(\lambda)_\mu$. Note that we have already seen that $(i,\mu)≥0$, and that $b_{i,l}z=0$ if $(i,\mu)=0$. Hence we assume $(i,\mu)>0$, otherwise there is nothing to prove. If $\omega_i≥2$, we have:\[
(\tilde f_{i,l}u,v)&=(b_{i,l}u,v)\\
&=-(u,K_{-li}a_{i,l}v)\\
&=-(u,K_{-li}a_{i,l}b_{i,\text{\sffamily c}}z)\\
&=-\bigg(u,K_{-li}\tau_{i,l}\sum_{k:\textup{\sffamily c}_k=l}\big(v_i^{-2l|\textup{\sffamily c}_{<k}|}K_{-li}-v_i^{2l|\textup{\sffamily c}_{<k}|}K_{li}\big)b_{i,\textup{\sffamily c}\backslash\textup{\sffamily c}_{k}}z\bigg)\\
&\qquad\qquad[\textit{c.f.~}\text{proof of~\ref{kashioplambda}}]\\
&\equiv-\bigg(u,\sum_{k:\textup{\sffamily c}_k=l}\big(v_i^{-2l|\textup{\sffamily c}_{<k}|}K_{-2li}-v_i^{2l|\textup{\sffamily c}_{<k}|}\big)b_{i,\textup{\sffamily c}\backslash\textup{\sffamily c}_{k}}z\bigg)\\
&\qquad\qquad[\textit{c.f.~}\ref{congrutau}]\\
&=-\bigg(u,\sum_{k:\textup{\sffamily c}_k=l}\big(v_i^{-2l|\textup{\sffamily c}_{<k}|}v^{-2l(i,\mu-|\textup{\sffamily c}\backslash\textup{\sffamily c}_{k}|i)}-v_i^{2l|\textup{\sffamily c}_{<k}|}\big)b_{i,\textup{\sffamily c}\backslash\textup{\sffamily c}_{k}}z\bigg)\\
&=-\bigg(u,\sum_{k:\textup{\sffamily c}_k=l}\big(v_i^{2l|\textup{\sffamily c}_{<k}|+4l|\textup{\sffamily c}_{>k}|}v^{-2l(i,\mu)}-v_i^{2l|\textup{\sffamily c}_{<k}|}\big)b_{i,\textup{\sffamily c}\backslash\textup{\sffamily c}_{k}}z\bigg)\\
&\equiv(u,b_{i,\textup{\sffamily c}\backslash\textup{\sffamily c}_{1}}z)\\
&=(u,\tilde e_{i,l}v).
\]
Thanks to the proof of~\ref{kashioplambda} and~\ref{congrutau}, the same can be proved if $\omega_i=1$. To that end, consider $v=b_{i,\nu}z$ and $u=b_{i,\nu'}z$ for some partitions $\nu,\nu'$ and elements $z,z'\in\mathcal K_i$, assuming again that $z\in V(\lambda)_\mu$. We have:\[
(\tilde f_{i,l}u,v)&=\sqrt\frac{l}{m_l(\nu')+1}(b_{i,l}u,v)\\
&=-\sqrt\frac{l}{m_l(\nu')+1}(u,K_{-li}a_{i,l}v)\\
&=-\sqrt\frac{l}{m_l(\nu')+1}(u,K_{-li}a_{i,l}b_{i,\nu}z)\\
&=-\sqrt\frac{l}{m_l(\nu')+1}\Big(u,K_{-li}\tau_{i,l}m_l(\nu)(K_{-li}-K_{li})b_{i,\nu\backslash l}z\Big)\\
&\qquad\qquad[\textit{c.f.~}\text{proof of~\ref{kashioplambda}}]\\
&\equiv-\sqrt\frac{l}{m_l(\nu')+1}\Big(u,\frac{m_l(\nu)}{l}(K_{-2li}-1)b_{i,\nu\backslash l}z\Big)\\
&\qquad\qquad[\textit{c.f.~}\ref{congrutau}]\\
&=-\sqrt\frac{l}{m_l(\nu')+1}\Big(u,\frac{m_l(\nu)}{l}(v^{-2l(i,\mu	)}-1)b_{i,\nu\backslash l}z\Big)\\
&\equiv\sqrt\frac{l}{m_l(\nu')+1}\frac{m_l(\nu)}{l}(u,b_{i,\nu\backslash l}z).\]
Iterating this computation, we see that $(\tilde f_{i,l}u,v)=(u,\tilde e_{i,l}v)=0$ if $\nu'\cup l\neq\nu$. Otherwise, we get:\[
(\tilde f_{i,l}u,v)&\equiv\sqrt\frac{m_l(\nu)}{l}(u,b_{i,\nu\backslash l}z)=(u,\tilde e_{i,l}v)
\]
 The case $i\in I^\text{re}$ is already known, but we reproduce the proof adapted to our conventions. The following can be proved by induction:\[
[E_i,F_i^{(n)}]=\tau_i\big(v^{-n+1}K_{-i}-v^{n-1}K_i)F_i^{(n-1)}.\]
Then, note that if $u=f_i^{m}u_0$ and $u'=f_i^{n}u'_0$, where $u_0,u'_0\in\mathcal K_i$, it is easy to prove that:\[
(\tilde f_i u,u')=(F_i^{(m+1)}u_0,F_i^{(n)}u'_0)=0\]
if $m+1\neq n$. If $n=m+1$ and $u'\in V(\lambda)_\mu$, we get:\[
(\tilde f_iu,u')&=\dfrac{1}{[m+1]}\big(F_iu,F_i^{(n)}u'_0\big)\\
&=\dfrac{1}{[n]}\Big(u,\dfrac{v}{v^{-1}-v}K_{-i}E_iF_i^{(n)}u'_0\Big)\\
&=\dfrac{1}{[n]}\Big(u,\dfrac{v}{v^{-1}-v}K_{-i}\tau_i\big(v^{-n+1}K_{-i}-v^{n-1}K_i)F_i^{(n-1)}u'_0\Big)\\
&\equiv\dfrac{1}{[n]}\Big(u,\dfrac{v}{v^{-1}-v}K_{-i}\big(v^{-n+1}K_{-i}-v^{n-1}K_i)F_i^{(n-1)}u'_0\Big)\\
&=\dfrac{1}{[n]}\Big(u,\dfrac{v}{v^{-1}-v}\big(v^{-n+1}v^{-2(i,\mu+i)}-v^{n-1})F_i^{(n-1)}u'_0\Big)\\
&=\Big(u,\dfrac{v^{-n+2}v^{-2(i,\mu+i)}-v^{n}}{v^{-n}-v^n}F_i^{(n-1)}u'_0\Big)\\
&=\Big(u,\dfrac{1-v^{-2n+2}v^{-2(i,\mu+i)}}{1-v^{-n}}F_i^{(n-1)}u'_0\Big)\\
&\equiv\big(u,F_i^{(n-1)}u'_0\big)\\
&=\big(u,\tilde e_iu'\big).
\]
We have assumed $(i,\mu+i)>0$, since otherwise we would have $u'=0$.
\edemo

The previous lemmas make it possible to reproduce step by step the original Kashiwara's grand loop argument (see~\cite[§4]{kashiloop}).

We also want to prove~\ref{jotens}, using the same kind of argument as in~\cite[4.7]{nakatens}. To that end, the characterization of the crystals $\mathcal B(\lambda)$ given by Joseph in~\cite[6.4.21]{Jo} has to be generalized. We first need two definitions:

\defi A crystal $\mathcal B$ is said to be of \emph{highest weight} $\lambda$ if:\benu
\item there exists $b_\lambda\in\mathcal B$ with weight $\lambda$ such that $\tilde e_\iota b_\lambda=0$ for every $\iota\in I_\infty$;
\item any element of $\mathcal B$ can be written $\tilde f_{\iota_1}\dots\tilde f_{\iota_r}b_\lambda$ for some $\iota_k\in I_\infty$.\eenu\edefi

\defi Consider a family $\{\mathcal B_\lambda\mid\lambda\in P^+\}$ of highest weight normal crystals $\mathcal B_\lambda$ of highest weight $\lambda$, with elements $b_\lambda\in\mathcal B_\lambda$ satisfying the properties of the above definition. It is called \emph{closed} if the subcrystal of $\mathcal B_\lambda\otimes \mathcal B_\mu$ generated by $b_\lambda\otimes b_\mu$ (\textit{i.e.}~obtained after successive applications of the $\tilde e_\iota$ and $\tilde f_\iota$) is isomorphic to $\mathcal B_{\lambda+\mu}$.\edefi
 
Our previous lemmas, along with the definitions and properties given in the previous sections make it possible to generalize the proof of~\cite[6.4.21]{Jo} and get:
 
 \prop The only closed family of highest weight normal crystals is $\{\mathcal B(\lambda)\mid\lambda\in P^+\}$.\eprop

Then it is easy to see that $\{\Irr\mathcal L(\lambda)\mid\lambda\in P^+\}$ is a closed family of highest weight normal crystals: thanks to~\ref{bijtens},~\ref{remnak} adapted to $\widetilde{\mathfrak Z}$ and~\ref{tenscrysgeom}, the arguments given in~\cite{nakatens} can be reproduced and we get~\ref{jotens}.
Alternatively (but similarly), the original proof given by Saito in~\cite{saito} can also be generalized to our framework.

\thanks{
\noindent 
\\
Department of Mathematics,
Massachusetts Institute of Technology
Cambridge, MA 02139, USA
\\ e-mail:\;\texttt{tbozec@mit.edu}}

\end{document}